\documentclass{article}

 \usepackage{graphicx}
  \usepackage{hyperref}
 \usepackage{amssymb,epsfig,amsmath,color,caption,subcaption,epstopdf,amsthm}
\usepackage{setspace}

\usepackage[all]{xy}

\newtheorem{definition}{Definition}[section]
\newtheorem{lemma}[definition]{Lemma}
\newtheorem{remark}[definition]{Remark}

\newtheorem{theorem}[definition]{Theorem}
\newtheorem{corollary}[definition]{Corollary}

\newcommand{\sctn}{\Gamma}

\begin{document}

\title{Feedback Integrators for \\ Nonholonomic Mechanical Systems}

\author{Dong Eui Chang and Matthew Perlmutter}

\date{June 30, 2018}

\maketitle

\begin{abstract}
The theory of feedback integrators is extended to handle mechanical systems with nonholonomic constraints with or without symmetry, so as to produce numerical integrators that preserve the nonholonomic constraints as well as other conserved quantities. To extend the feedback integrators, we develop a suitable extension theory for nonholonomic systems, and also a corresponding reduction theory for systems with symmetry. It is then applied to various nonholonomic systems such as 
the Suslov problem on $\operatorname{SO}(3)$, the knife edge, the Chaplygin sleigh, the vertical rolling disk, the roller racer,  the Heisenberg system, and the nonholonomic oscillator. 
\end{abstract}

\tableofcontents

\section{Introduction}

The feedback integrators method was  rigorously developed in \cite{ChJiPe16} to numerically integrate the equations of motion of dynamical systems so that the conserved quantities or first integrals of a given system are numerically well preserved. The paradigm of the method is not to develop any specific discrete integration scheme, but to modify a given continuous-time system so that any off-the-shelf integrator, such as Euler or Runge-Kutta, can be applied to integrate the equations of motion of the modified system. Its excellent performance was well demonstrated in \cite{ChJiPe16} on the three benchmark examples: the free rigid body system, the two-body problem and a perturbed two-body problem, all of which are  systems with no constraints. In this paper, we extend the theory of feedback integrators to mechanical systems with nonholonomic constraints for which we take the Hamiltonian approach rather than the Lagrangian approach, since we can express the Lagrange multiplier for the constraint in a simple form using the Poisson bracket.

We note that nonholonomic systems are important in mechanics, robotics, and control, and as their dynamical and geometric properties don't satisfy the symplecticity of Hamiltonian systems, they continue to be studied from various points of view. An important reference is \cite{Bl03}. We note here also several papers, including some recent, on the geometry and dynamics of nonholonomic systems,
\cite{Ko, vaMa, CaLeDi, GaNa}.



There is by now a significant literature on the problem of extending symplectic integrators to geometric (variational) integrators on the class of nonholonomic systems.  The discrete Lagrange d'Alembert equations (DLA equations), extending the discrete variational principle, and corresponding discrete Euler Lagrange equations (\cite{MaWe}) to the case of nonholonomic velocity constraints, first appeared in the paper of \cite{Co}. The theory was then further developed in \cite{McPe, FeZen, IgMarMart, KoMaSu}. These papers also consider reduction of the corresponding discrete flow and variational principle under certain natural settings, namely the Suslov problem (Lie group configuration space), and the Chaplygin systems where the symmetry group orbits are transverse to the velocity constraints.
More recently, the theory of DLA equations has been further developed in the works of \cite{FeIgMa08,FeJiMa,JiSc} where the GNI (geometric nonholonomic integrator), derived with the added structure of a metric and its corresponding projectors on the nonholonomic constraint distribution, is derived and its convergence properties are studied.


By contrast to this approach, we will pursue the notion of feedback integrators and adapt them to nonholonomic systems. A fundamental difference is that we work on the Hamiltonian side, first developing an extension theory for nonholonomic systems with the advantage that once we extend the vector field to all of $T^{\ast}Q$ we can apply the methodology of feedback integrators to the extended system.

This theory of nonholonomic feedback integrators can be summarized as follows. Given a mechanical system $\Sigma$ on a cotangent bundle $T^*Q$ with a set $\mathcal C$ of nonholonomic constraints, first extend the system $\Sigma$ from $\mathcal C$  to its full phase space $T^*Q$  such that the extended system, denoted $\tilde \Sigma$, is free of any nonholonomic constraint but has $\mathcal C$ as its invariant set. Since the constraint $\mathcal C$ has been eliminated, the ordinary feedback integrators, which were developed in \cite{ChJiPe16}, can be applied to $\tilde \Sigma$. In other words, embed $\tilde \Sigma$ from $T^*Q$ into some $\mathbb R^n$ if $T^*Q$ is not homeomorphic to Euclidean space, and then extend the system $\tilde \Sigma$ to a neighbourhood of $T^*Q$ in $\mathbb R^n$, where the extended system is denoted $\tilde \Sigma_e$.  Then, collect conserved quantities of $\tilde \Sigma_e$ including the nonholonomic constraint set $\mathcal C$ and pick  a point $x_0$ from $\mathcal C \subset T^*Q \subset \mathbb R^n$   which is an initial point of a trajectory to be computed.   
Modify the dynamics of $\tilde \Sigma_e$ outside the intersection, denoted $\Lambda$, of the level sets of the conserved quantities containing $x_0$, while leaving the dynamics on $\Lambda$ unchanged, such that $\Lambda$ becomes an attractor of the modified dynamics. Finally, integrate the resultant dynamics from the point $x_0$ using any off-the-shelf integrator such as Euler, Runge-Kutta, etc. This integrator will then compute a numerical trajectory that stays close to the set $\Lambda$, thus  preserving all the conserved quantities. The trajectory converges to $\Lambda$ as the size of integration time step tends to zero \cite{ChJiPe16}. Two main advantages of  feedback integrators are the use of one single global Cartesian coordinate system in the ambient space $\mathbb R^n$ and the use of any off-the-shelter numerical integrator.

Furthermore, this extension theory of nonholonomic systems from $\mathcal{C}$ to $T^{\ast}Q$ respects symmetries: if the nonholonomic constraint is invariant under a Lie group action, then the extended Hamiltonian and the extended vector field are also shown to be invariant and therefore induces a vector field on the quotient. We will also show in the case of nonholonomic systems on Lie groups and on trivial principal bundles, that, in the presence of symmetry,  extension commutes with reduction: namely the extended Hamiltonian and vector field drops to the quotient and is equal to the extension of the reduced nonholonomic and Hamiltonian data.  In an upcoming paper we will show that these results hold generally on $T^{\ast}Q/G$, although the presence of curvature makes the reduction theory more involved. 

A key consequence of this reduction theory is that our feedback integrators for nonholonomic systems naturally drop to the quotient space, and the reduced integrator can be directly constructed from the reduced extension vector field. We will see this procedure in the examples that are considered in the second half of the paper, in particular for the Chaplygin sleign and the roller racer.

In the special case of a nonholonomic system on a Lie group, the reduction and reconstruction equations of the nonholonomic system leads to a beautiful generalization of the Lie Poisson bracket. The reduced vector field consists of the sum of a Lie Poisson extended Hamiltonian vector field 
and the constraint forces appearing as a linear combination of the covectors that determine the nonholonomic constraint.
When these vanish, one recovers the Lie Poisson equations.


The Theorems proven in sections 2.2 and 2.3, that the process of extension commutes with reduction, immediately has an application to the feedback integrators for nonholonomic systems with symmetry. Namely, the feedback integrator for the unreduced extended system will project under the quotient to the feedback integrator for the reduced system.


This paper is organized as follows. First, the theory of extending nonholonomic systems from their constraints to  their ambient space is developed for the canonical case, i.e. when the phase space is a cotangent bundle equipped with the canonical bracket. We then consider two extreme cases for symmetric systems. The first one is the case where  the configuration space $Q$ is the symmetry group of  system, and the second is where the configuration space is a trivial principal bundle with its fiber as the symmetry group, leaving the nontrivial bundle case as future work. Next, the theory of feedback integrators is briefly reviewed from \cite{ChJiPe16}, and then finally the nonholonomic feedback integrators are constructed for the following systems: the Suslov problem on $\operatorname{SO}(3)$, the knife edge, the Chaplygin sleigh, the vertical rolling disk, the roller racer,  the Heisenberg system, and the nonholonomic oscillator. For the purpose of evaluation, we compare the trajectories generated by feedback integrators with exact solutions and the trajectories generated by other integrators for some of the systems, demonstrating the efficacy of feedback integrators.

\section{Extension of Nonholonomic Mechanical Systems}

In this section, we develop a theory of extending nonholonomic mechanical systems from their constraints to the ambient space, so as to effectively design feedback integrators.    After the development of the extension theory, we explain how to design feedback integrators for numerical integration of the extended system with preservation of conserved quantities including the Hamiltonian and the constraint set.


\subsection{Noholonomic Systems on Cotangent Bundles}
Let the phase space be the cotangent bundle $T^*Q$ of a configuration space $Q$, and use coordinates $(q,p)$ for $T^*Q$.  Consider a Hamiltonian
\begin{equation}\label{def:H}
H(q,p) = \frac{1}{2}\langle p,  m(q)^{-1}p\rangle + V(q) = \frac{1}{2}\langle p, p^\sharp\rangle  + V(q),
\end{equation}
where $m$ is the symmetric positive definite mass tensor on $Q$ and $V$ is the potential energy of the Hamiltonian.
The musical maps $\sharp: T^*Q \rightarrow TQ$ and $\flat :  TQ\rightarrow T^*Q$ are understood with respect to  the metric $m$ on $Q$.  For example, $p^\sharp = m^{-1}p$ and $v^\flat = mv$ for $p \in T^*Q$ and $v\in TQ$. Consider a nonholonomic constraint set  
\begin{equation}\label{non:holo:set}
\mathcal C =\{ (q,p) \in T^*Q \mid\langle p, e_i(q) \rangle = 0, i = 1, ..., K\},
\end{equation}
 where $\{e_1, \ldots, e_K\}$ is a set of  orthonormal vector fields on $Q$ with respect to the metric $m$.

The equations of motion of the mechanical system with the Hamiltonian $H$ and  the constraint $\mathcal C$ are given by
\begin{subequations}\label{original:nonholo}
\begin{align}
\dot q &= \frac{\partial H}{\partial p}, \label{bam:q}\\
\dot p &= -\frac{\partial H}{\partial q} + \sum_{i=1}^K  \lambda_i e_i^\flat,\label{bam:p}\\
(q,p) &\in \mathcal C, 
\end{align}
\end{subequations}
where the multipliers $ \lambda_i$'s are determined to make each $\langle p, e(q)_i\rangle$ a constant of motion  of \eqref{bam:q} and \eqref{bam:p} on $T^*Q$, i.e., 
\[
\frac{d}{dt}\langle p(t), e_i(q(t)) \rangle = 0
\]
for all $t\in \mathbb R$ and  all $i = 1,\ldots, K$ along the flow of \eqref{bam:q} and \eqref{bam:p}. It follows that
\begin{align}\label{def:lambda}
\lambda_i (q,p)&=   \langle \frac{\partial H}{\partial q},e_i\rangle - \langle p, \mathrm{D}e_i(q) \cdot \frac{\partial H}{\partial p}\rangle = \{ H, \langle p, e_i \rangle\}
\end{align}
for $i = 1, \ldots, K$. 
This choice of multipliers makes $\mathcal C$ invariant under the flow of the system \eqref{bam:q} and \eqref{bam:p}.
One can verify that
\[
\frac{dH}{dt}= \sum_{i=1}^K\lambda_i \langle p, e_i\rangle,
\]
which vanishes on $\mathcal C$, but  not necessarily everywhere on $T^*Q$. However, we  want a Hamiltonian function to be a constant of motion in the entire phase space. 

Let us extend the system from $\mathcal C$ to the entire phase space $T^*Q$. First, define an ``extended'' Hamiltonian function $\tilde H$ by
\begin{align}\label{def:H:tilde}
\tilde H &= \frac{1}{2}\left \langle (p -\sum_{j=1}^K \langle p, e_j\rangle e_j^\flat),  (p -\sum_{j=1}^K \langle p, e_j \rangle e_j^\flat )^\sharp \right \rangle + V(q) \nonumber \\
&=  \frac{1}{2}\langle p, p^\sharp\rangle   - \frac{1}{2}\sum_{j=1}^K |\langle p, e_j\rangle|^2 + V(q),
\end{align}
which coincides with the original Hamiltonian $H$ on the nonholonomic constraint set $\mathcal C$.

Consider the following {\it extended} system on $T^*Q$:
\begin{subequations}\label{bam:qp}
\begin{align}
\dot q &= \frac{\partial \tilde H}{\partial p} =  p^\sharp - \sum_{j=1}^K \langle p, e_j\rangle e_j,\label{bambam:q}\\
\dot p &= - \frac{\partial \tilde H}{\partial q} + \sum_{i=1}^K\tilde \lambda_i (q,p)e_i^\flat \label{bambam:p}
\end{align}
\end{subequations}
 where
\begin{align}\label{def:lambda:tilde}
\tilde \lambda_i (q,p)&= \{ \tilde H,  \langle p, e_i \rangle\}\\
&=  \{ H,\langle p, e_i \rangle \} - \sum_{j=1}^K   \langle p, e_j \rangle \{ \langle p, e_j \rangle, \langle p, e_i \rangle\} \nonumber\\
 &=  \{ H,\langle p, e_i \rangle \} + \sum_{j=1}^K   \langle p, e_j \rangle  \langle p, [e_j,e_i] \rangle \rangle.\nonumber
\end{align}

\begin{theorem}\label{theorem:extension:nonhol}
The extended system \eqref{bam:qp} coincides on the constraint set $\mathcal C$ with the constrained system \eqref{original:nonholo}. Moreover, the  constraint set $\mathcal C$ is invariant under the flow of \eqref{bam:qp}. 
\begin{proof}
Straightforward.
\end{proof}
\end{theorem}

\begin{theorem}\label{theorem:first:integrals}
The Hamiltonian $\tilde H$  given in \eqref{def:H:tilde} and the constraint momentum components $\langle p, e_i\rangle$, $i=1, \ldots,K$ are  the first integrals of the extended system  \eqref{bam:qp} on  $T^*Q$. Moreover, $\tilde H = H$ and $\tilde \lambda_i = \lambda_i$, $i=1, \ldots, K$ on $\mathcal C$, where $H$ is given in (\ref{def:H}) and  $\tilde \lambda_i$ and $\lambda_i$ are given in \eqref{def:lambda:tilde} and \eqref{def:lambda}.
\begin{proof}
Along the flow of   \eqref{bam:qp}
\begin{align*}
\frac{d \tilde H}{dt} &= \frac{\partial \tilde H}{\partial q} \dot q + \frac{\partial \tilde  H}{\partial p}\dot p\\
&=  \frac{\partial \tilde H}{\partial q}\frac{\partial \tilde H}{\partial p} + \frac{\partial \tilde H}{\partial p}(-\frac{\partial \tilde H}{\partial q} + \sum_{i=1}^K \tilde \lambda_i e_i^\flat)\\
&= -   \sum_{i=1}^K \tilde \lambda_i e_i^\flat \frac{\partial \tilde H}{\partial p}\\
&= -   \sum_{i=1}^K \tilde \lambda_i \langle e_i^\flat, p^\sharp - \sum_{j=1}^K \langle p, e_j\rangle e_j\rangle \\
&=0
\end{align*}
since $\{e_i\}$ is a set of orthonormal vector fields with respect to the metric $m$. Hence, $\tilde H$ is a first integral of \eqref{bam:qp} on $T^*Q$. The other statements are straightforward to prove.

\end{proof}
\end{theorem}

\begin{corollary}\label{corollary:H:instead}
The original Hamiltonian $H$  in \eqref{def:H}  can be written as
\begin{equation}\label{H:intermsof:Htilde:pe}
H = \tilde H + \frac{1}{2}\sum_{i=1}^K|\langle p, e_i\rangle|^2
\end{equation}
and it is a first integral of the extended system \eqref{bam:qp} on $T^*Q$.
\begin{proof}
Equation \eqref{H:intermsof:Htilde:pe} follows from \eqref{def:H} and \eqref{def:H:tilde}. 
Since $\tilde H$ and $\langle p, e_i \rangle$, $i = 1, \ldots, K$ are the first integrals of   \eqref{bam:qp},  $H$ is also a first integral of  \eqref{bam:qp} by \eqref{H:intermsof:Htilde:pe}.
\end{proof}
\end{corollary}
We emphasize that  the system \eqref{bam:qp} is not subject to the constraint $\mathcal C$ any more, but  is an ordinary  Hamiltonian system on $T^*Q$ with the force $ \sum_{i=1}^K\tilde \lambda_i (q,p)e_i^\flat$ that is gyroscopic, i.e. not affecting the value of the Hamiltonian $\tilde H$. 

\begin{remark}  \label{remark:nonhol:general}
Let us consider a more general case where the constraint vector fields $e_i$'s are not necessarily orthonormal with respect to the kinetic energy metric.  In such a case, the multipliers $\lambda_i$, $i = 1, \ldots, K$, are computed as 
\begin{equation}\label{lambda:general}
\lambda_i = \sum_{\ell = 1}^K\{H, \langle p, e_\ell\rangle \}C_{\ell i},
\end{equation}
where $C_{\ell i}$ is the inverse matrix of the symmetric matrix
\[
C^{\ell i} = \langle e_\ell^b, e_i\rangle
\]
 so as to make  each constraint momentum $\langle p, e_i\rangle$ for $i = 1, \ldots, K$,  a first integral of \eqref{bam:q} and \eqref{bam:p} on $T^*Q$. The extended Hamiltonian $\tilde H$ is modified as
\begin{align}\label{def:H:tilde:general}
\tilde H (q,p) &=  \frac{1}{2}\left \langle (p -\sum_{i=1}^K \sum_{j = 1}^K \langle p, e_j \rangle C_{ ji}e_i^\flat),   (p - \sum_{i=1}^K \sum_{j = 1}^K \langle p, e_j \rangle C_{ ji}e_i^\flat)^\sharp  \right \rangle + V(q)  \nonumber\\
&= H(q,p) - \frac{1}{2} \sum_{i=1}^K \sum_{j=1}^K \langle p,e_i \rangle C_{ij}\langle p, e_j\rangle,
\end{align}
The multipliers $\tilde \lambda_i$'s are modified as
\begin{equation}\label{lambda:tilde:general}
\tilde \lambda_i = \sum_{\ell = 1}^K\{\tilde H, \langle p, e_\ell\rangle \}C_{\ell i}.
\end{equation}
Equation \eqref{bambam:q} is modified accordingly as
\[
\dot q = \frac{\partial \tilde H}{\partial q} = p^\sharp - \sum_{j=1}^K\sum_{\ell=1}^K \langle p, e_j \rangle C_{j \ell} e_{\ell},
\]
and equation \eqref{H:intermsof:Htilde:pe} is modified to
\[
H = \tilde H +  \frac{1}{2} \sum_{i=1}^K \sum_{j=1}^K \langle p,e_i \rangle C_{ij}\langle p, e_j\rangle.
\]
Then,  all the theorems and corollaries in the above  hold true in this more general form. The verification is left to the reader. We however do not use these general formulas in this paper.
\end{remark}

\subsection{Nonholonomic Systems on Lie Groups}
We consider a nonholonomic system on the Lie group $G$. We shall assume that the noholonomic constraint distribution, $\mathcal{C}$ is left $G$-invariant and is therefore determined by a subspace $\mathcal{C}_{\rm red}\subset \mathfrak{g}^*$. Using concatenated notation for the tangent/cotangent lift of the group multiplication, we have that
\begin{equation}
\label{distribution}
\mathcal{C}(g)=g\cdot \mathcal{C}_{\rm red}\subset T_{g}^*G.
\end{equation}
Suppose that we are given a  $G$-invariant Hamiltonian $H\in C^\infty(T^*G)$, which induces its corresponding reduced Hamiltonian $h \in C^\infty(\mathfrak{g}^*)$ such that $h\circ \pi = H$ where $\pi : T^*G \rightarrow \mathfrak{g}^*$ just the group projection for the left action, given by $\pi (\alpha_g) = g^{-1}\cdot \alpha_g$. The reduced Hamiltonian can then be  written as 
\begin{equation}\label{def:reduced:h}
h(\mu) = \frac{1}{2}\langle \mu, {\mathbb I}^{-1}\mu \rangle =  \frac{1}{2}\langle \mu, \mu^\sharp \rangle
\end{equation}
where $\mathbb I: \mathfrak g \rightarrow \mathfrak g^*$ is the locked inertia tensor  which is symmetric and positive definite.  We here assume that $H$ has kinetic energy only. The musical maps $\sharp: \mathfrak g^* \rightarrow \mathfrak g$ and $\flat : \mathfrak g \rightarrow \mathfrak g^*$ are understood with respect to  $\mathbb I$. For example, $\mu^\sharp = \mathbb I^{-1}\mu$ and $e^\flat = \mathbb I e$ for all $\mu \in\mathfrak g^*$ and $e\in \mathfrak g$.  With this locked inertia tensor,  the reduced nonholonomic constraint set can be written as
\[
{\mathcal C}_{\rm red} = \{ \mu \in {\mathfrak g}^* \mid  \langle \mu, e_i\rangle = 0, i = 1, \ldots, K\},
\]
where $\{e_1, \ldots, e_K\} \subset \mathfrak g$ is a set of orthonormal Lie algebra elements with respect to the metric $\mathbb I$, i.e. $\langle \mathbb Ie_i, e_i\rangle = 1$. Here, each $e_i$ can be understood as the valuation of the corresponding $e_i (g)$  in \eqref{non:holo:set} at the identity element of $G$.

By the theory of Lie-Poisson reduction \cite{MaRa02},  the set of equations of motion on $G\times \mathfrak{g}^*$ is given by
\begin{subequations}\label{sys:LiePoisson}
\begin{align}
\displaystyle
\dot g &=  g\cdot \frac{\delta h}{\delta \mu}= g\cdot {\mathbb I}^{-1}\mu, \label{sys:LiePoisson:g}\\
 \dot \mu &= \operatorname{ad}^*_{\frac{\delta h}{\delta \mu}} \mu + \sum_{i=1}^K  \lambda_i  e_i^b ,\label{sys:LiePoisson:mu}\\
\mu &\in {\mathcal C}_{\rm red},
\end{align}
\end{subequations}
where the final two equations are the reduced nonholonomic system on $\mathcal{C}_{\rm red}\subset \mathfrak{g}^{\ast}$. They are solved for $\mu(t)$, which then determines the curve $\frac{\delta h}{\delta \mu}(t)\in \mathfrak{g}$, which in turn determines the differential equation on $G$, \eqref{sys:LiePoisson:g}.
The multipliers, $\lambda_i$'s, are determined by the condition that each $\langle \mu, e_i\rangle$ for $i= 1, \ldots, K$  is a first integral of \eqref{sys:LiePoisson:mu} so as to make $\mathcal C_{\rm red}$ an invariant set of the system. It follows that
\begin{equation}\label{red:lambda}
\lambda_i(\mu) =- \langle \operatorname{ad}^*_{\frac{\delta h}{\delta \mu}} \mu, e_i \rangle =  -\left \langle \mu, \left  [ \frac{\delta h}{\delta \mu}, e_i \right ] \right \rangle = \{h, e_i\}_{-}(\mu),
\end{equation}
where $[ \, , ]$ is the bracket on the Lie algebra $\mathfrak g$ and 
 $\{\, ,\}_{-}$ is  the minus Lie-Poisson bracket on $\mathfrak g^*$ that, recall, is given by
\[
\{f_1,f_2\} (\mu) = - \left \langle \mu,  \left [ \frac{\delta f_1}{\delta \mu}, \frac{\delta f_2}{\delta \mu}\right ]\right \rangle
\]
for $f_1, f_2\in C^\infty({\mathfrak g}^*)$. Notice that  the expression for $\lambda_i$ in \eqref{red:lambda} is exactly the reduced version of \eqref{def:lambda} as it should be.

Let us now extend this system from $G\times \mathcal C_{\rm red}$ to the entire phase space $T^*G  = G \times \mathfrak g^*$.
First, we define an extended Hamiltonian $\tilde h : \mathfrak g^* \rightarrow \mathbb R$ by
\begin{align}\label{def:reduced:h:tilde}
\tilde h(\mu) &=  \frac{1}{2}\langle (\mu - \sum_{i=1}^K \langle \mu, e_i\rangle e_i^\flat ) , {\mathbb I}^{-1} (\mu - \sum_{i=1}^K \langle \mu, e_i\rangle e_i^\flat )  \rangle \nonumber \\
&=  h(\mu) - \frac{1}{2} \sum_{i=1}^K |\langle \mu, e_i\rangle |^2.
\end{align}
Consider the following new dynamical system on $G \times \mathfrak g^*$:
\begin{subequations}\label{new:red:bam}
\begin{align}
\dot g &=  g\cdot \frac{\delta \tilde h}{\delta \mu}= g\cdot ( {\mathbb I}^{-1}\mu - \sum_{i=1}^K \langle \mu, e_i\rangle e_i),\label{new:red:g}\\
\dot \mu &= \operatorname{ad}^*_{\frac{\delta \tilde h}{\delta \mu}} \mu + \sum_{i=1}^K \tilde \lambda_i  e_i^b ,\label{new:red:mu}
\end{align}
\end{subequations}
where
\begin{equation} \label{def:reduced:lambda:tilde}
\tilde \lambda_i (\mu) =  \{\tilde h, \langle \mu, e_i\rangle\}_{-} =\{ h, \langle \mu, e_i\rangle\}_{-} + \sum_{j=1}^K \langle \mu, e_j\rangle \langle \mu, [e_j, e_i]\rangle.
\end{equation}
Notice that $\tilde \lambda_i$ in \eqref{def:reduced:lambda:tilde} is the reduced version of \eqref{def:lambda:tilde} and that the system \eqref{new:red:bam} is the left-trivialization or the Lie-Poisson reduction of \eqref{bam:qp}. We have the following commutative diagram asserting that extension commutes with reduction:
\[
\xymatrixcolsep{4pc}
\xymatrix{
(H, \{\lambda_i\}, \mathcal C) \ar[d]_{\textup{reduction}} \ar[r]^{\textup{extension}} & (\tilde H, \{\tilde \lambda_i\},G\times \mathfrak g^*) \ar[d]^{\textup{reduction}} \\
(h, \{\lambda_i\},{\mathcal C}_{\rm red}) \ar[r]^{\textup{extension}} & (\tilde h, \{\tilde \lambda_i\}, \mathfrak g^*)}
\]
where $\{\lambda_i\}$ and $\{\tilde \lambda_i\}$ represent the corresponding forces. 
This commutative diagram, together with Theorems \ref{theorem:extension:nonhol} and \ref{theorem:first:integrals}, and Corollary \ref{corollary:H:instead},  implies the following results:
\begin{theorem}\label{theorem:reduced:extension:nonhol}
The extended system \eqref{new:red:bam} coincides on the nonholonomic constraint set  with the constrained system \eqref{sys:LiePoisson}. Moreover, the  constraint set $\mathcal C_{\rm red}$ is invariant under the flow of \eqref{new:red:bam}. 
\end{theorem}

\begin{theorem}\label{theorem:reduced:first:integrals}
The reduced Hamiltonian $\tilde h$  in \eqref{def:reduced:h:tilde} and the constraint momentum components $\langle \mu, e_i\rangle$, $i=1, \ldots,K$ are the first integrals of the extended system \eqref{new:red:bam} on  $G \times \mathfrak g^*$. Moreover, $\tilde h = h$ and $\tilde \lambda_i = \lambda_i$, $i=1, \ldots, K$ on $\mathcal C_{\rm red}$, where $h$ is given in \eqref{def:reduced:h} and  $\tilde \lambda_i$ and $\lambda_i$ are given in \eqref{def:reduced:lambda:tilde} and \eqref{red:lambda}.
\end{theorem}

\begin{corollary}\label{corollary:reduced:H:instead}
The original reduced Hamiltonian $h$  in \eqref{def:reduced:h}  can be written as
\[
h = \tilde h + \frac{1}{2}\sum_{i=1}^K|\langle \mu, e_i\rangle|^2
\]
and it is a first integral of the extended system \eqref{new:red:bam} on  $G \times \mathfrak g^*$.
\end{corollary}

\subsection{Nonholonomic Mechanical Systems with Symmetry on Trivial Principal Bundles}

Consider the case in which the configuration space $Q$  is the product of a Lie group $G$ and a manifold $X$, i.e. $Q = G \times X$.  The group $G$ acts on $Q$ by left multiplication on the first factor $G$ of $Q$. By the cotangent lifted action of $G$, we have $T^*Q \simeq T^*G \times T^*X \simeq G \times {\mathfrak g}^* \times T^*X$. Use $(g,\mu, x, p_x)$ for coordinates on $G \times {\mathfrak g}^* \times T^*X$.  According to  \cite{MoMaRa84}, the Poisson bracket on $ {\mathfrak g}^* \times T^*X$ induced from the canonical bracket on $T^*Q$ is given by
\begin{equation}\label{bracket:QG:trivial}
\{f_1, f_2\} (\mu,x,p_x) =  - \left \langle \mu,  \left [ \frac{\delta f_1}{\delta \mu}, \frac{\delta f_2}{\delta \mu}\right ]\right \rangle + \frac{\partial f_1}{\partial x}\cdot \frac{\partial f_2}{\partial p_x}  - \frac{\partial f_2}{\partial x}\cdot \frac{\partial f_1}{\partial p_x}
\end{equation}
for all $f_1, f_2\in C^\infty({\mathfrak g}^* \times T^*X)$, where $[\, ,]$ is the bracket on the Lie algebra $\mathfrak g$. Sometimes, the bracket in \eqref{bracket:QG:trivial} is compactly written as
\[
\{f_1,f_2\}_{{\mathfrak g}^* \times T^*X} = \{f_1,f_2\}_{{\mathfrak g}^*} + \{f_1,f_2\}_{T^*X}. 
\]

Let $H \in C^\infty(T^*Q)$ be a  $G$-invariant Hamiltonian of the form \eqref{def:H}. It induces a reduced Hamiltonian $h \in C^\infty({\mathfrak  g}^* \times T^*X)$ that can be written as
\begin{equation}\label{def:reduced:h:QG:trivial}
h(\mu, x,p_x) = \frac{1}{2} \langle (\mu, p_x), m(x)^{-1}\cdot (\mu, p_x) \rangle + V(x),
\end{equation}
where $m$ is the reduced mass tensor and $V$ is the reduced potential energy of the system. 
Let $\mathcal C \in \sctn (T^*Q)$ be a nonholonomic constraint set that is $G$-invariant.  Then, it induces  a reduced nonholonomic constraint set ${\mathcal C}_{\rm red}$ which can be written as 
\[
{\mathcal C}_{\rm red} = \{ (\mu, p_x) \in \mathfrak g^* \times T^*X \mid  \langle (\mu, p_x), e_i\rangle = 0, i = 1,\ldots, K\}
\]
where $\{e_i\}$ is a set of sections of ${\mathfrak g}^* \times TX$ over $X$, or loosely speaking, vector fields on ${\mathfrak g}^* \times X$, that are orthonormal with respect to the reduced mass tensor $m$.

The equations of motion of the constrained Hamiltonian system with the Hamiltonian $h$ and the constraint ${\mathcal C}_{\rm red}$ are given by
\begin{subequations}\label{sys:QG:trivial}
\begin{align}
\displaystyle
(\dot g, \dot x) &=   \left (g\cdot \frac{\delta h}{\delta \mu},
\frac{\partial h}{\partial p_x}\right ), \label{sys:QG:trivial:q}\\
( \dot \mu,
 \dot p_x)
 &= 
\left (
 \operatorname{ad}^*_{\frac{\delta h}{\delta \mu}} \mu, 
\displaystyle -\frac{\partial h}{\partial x}
\right )
+ \sum_{i=1}^K  \lambda_i  e_i^b ,\label{sys:QG:trivial:p}\\
(\mu, p_x) &\in {\mathcal C}_{\rm red},
\end{align}
\end{subequations}
where 
\[
e_i^\flat = me_i
\]
for $i=1, \ldots, K$.
The multipliers $\lambda_i$'s, $i=1, \ldots, K$, are determined to make each constraint momentum $ \langle (\mu, p_x), e_i\rangle$ a first integral of the unconstrained system in \eqref{sys:QG:trivial:q} and \eqref{sys:QG:trivial:p} so as to make the constraint set ${\mathcal C}_{\rm red}$ an invariant set. Hence,
\begin{align}\label{def:lambda:QG:trivial}
\lambda_i &=   \{ h, \langle (\mu,p_x), e_i \rangle\}
\end{align}
for $i =1, \ldots, K$, where the bracket is the one given in \eqref{bracket:QG:trivial}. This is exactly the reduced version of \eqref{def:lambda}.

Let us now extend this system from $G\times \mathcal C_{\rm red}$ to the entire phase space $T^*Q  = G \times \mathfrak g^* \times T^*X$.
First,  define an extended Hamiltonian $\tilde h \in C^\infty(\mathfrak g^* \times T^*X)$ by
\begin{align}\label{def:reduced:h:tilde:QG:trivial}
\tilde h(\mu,x,p_x) 
&=  h(\mu,x,p_x) - \frac{1}{2} \sum_{i=1}^K |\langle (\mu,p_x), e_i\rangle |^2.
\end{align}
Consider the following new dynamical system on $T^*Q = G \times \mathfrak g^* \times T^*X$:
\begin{subequations}\label{new:red:bam:QG:trivial}
\begin{align}
\displaystyle
(\dot g, \dot x) &=   \left (g\cdot \frac{\delta \tilde h}{\delta \mu},
\frac{\partial \tilde h}{\partial p_x}\right ), \label{new:sys:QG:trivial:q}\\
( \dot \mu,
 \dot p_x)
 &= 
\left (
 \operatorname{ad}^*_{\frac{\delta \tilde h}{\delta \mu}} \mu, 
\displaystyle -\frac{\partial \tilde h}{\partial x}
\right )
+ \sum_{i=1}^K  \tilde \lambda_i  e_i^b ,\label{new:sys:QG:trivial:p}
\end{align}
\end{subequations}
where
\begin{align} \label{def:reduced:lambda:tilde:QG:trivial}
\tilde \lambda_i  &=  \{\tilde h, \langle(\mu,p_x), e_i\rangle\}. 
\end{align}
Notice that $\tilde \lambda_i$ in \eqref{def:reduced:lambda:tilde:QG:trivial} is the reduced version of \eqref{def:lambda:tilde} and that the system \eqref{new:red:bam:QG:trivial} is the left-trivialization  of \eqref{bam:qp}. We again have the following commutative diagram:
\[
\xymatrixcolsep{4pc}
\xymatrix{
(H, \{\lambda_i\}, \mathcal C) \ar[d]_{\textup{reduction}} \ar[r]^{\textup{extension\hspace{10mm}}} & (\tilde H, \{\tilde \lambda_i\},G\times \mathfrak g^* \times T^*X) \ar[d]^{\textup{reduction}} \\
(h, \{\lambda_i\},{\mathcal C}_{\rm red}) \ar[r]^{\textup{extension\hspace{5mm}}} & (\tilde h, \{\tilde \lambda_i\}, \mathfrak g^* \times T^*X)}
\]
where $\{\lambda_i\}$ and $\{\tilde \lambda_i\}$ represent the corresponding forces. 
This commutative diagram, together with Theorems \ref{theorem:extension:nonhol} and \ref{theorem:first:integrals}, and Corollary \ref{corollary:H:instead}, implies  the following results:
\begin{theorem}\label{theorem:reduced:extension:nonhol:QG:trivial}
The extended system \eqref{new:red:bam:QG:trivial} coincides on the nonholonomic constraint set  with the constrained system \eqref{sys:QG:trivial}. Moreover, the  constraint set $\mathcal C_{\rm red}$ is invariant under the flow of \eqref{new:red:bam:QG:trivial}.
\end{theorem}

\begin{theorem}\label{theorem:reduced:first:integrals:QG:trivial}
The reduced Hamiltonian $\tilde h$  in \eqref{def:reduced:h:tilde:QG:trivial} and the constraint momentum components $\langle ( \mu,p_x), e_i\rangle$, $i=1, \ldots,K$ are the first integrals of the extended system \eqref{new:red:bam:QG:trivial} on  $T^*Q$. Moreover, $\tilde h = h$ and $\tilde \lambda_i = \lambda_i$, $i=1, \ldots, K$ on $\mathcal C_{\rm red}$, where $h$ is given in \eqref{def:reduced:h:QG:trivial} and  $\tilde \lambda_i$ and $\lambda_i$ are given in \eqref{def:reduced:lambda:tilde:QG:trivial} and \eqref{def:lambda:QG:trivial}.
\end{theorem}

\begin{corollary}\label{corollary:reduced:H:instead:QG:trivial}
The original reduced Hamiltonian $h$  in \eqref{def:reduced:h:QG:trivial}  can be written as
\[
h = \tilde h + \frac{1}{2}\sum_{i=1}^K|\langle (\mu,p_x), e_i\rangle|^2
\]
and it is a first integral of the extended system \eqref{new:red:bam:QG:trivial} on  $T^*Q$.
\end{corollary}

\subsection{Design of Feedback Integrators for Nonholonomic Mechanical Systems}
We first briefly review the  theory of feedback integrators from \cite{ChJiPe16}; refer to \cite{ChJiPe16} for more detail.
Consider a dynamical system on an open subset $U$ of $\mathbb R^n$:
\begin{equation}\label{our:dyn:sys}
\dot x = X(x),
\end{equation}
where $X$ is a $C^1$ vector field on $U$.  We make the following three assumptions:
\begin{itemize}
\item [A1.]  There is a $C^2$ function $V : U \rightarrow \mathbb R$ such that  $V(x) \geq 0$ for all $x\in U$, $V^{-1}(0) \neq \emptyset$, and it is a first integral of \eqref{our:dyn:sys}, i.e.
\begin{equation}\label{Vdot:X}
\nabla V(x) \cdot X (x) = 0
\end{equation}
for all $x \in U$. 
 
 \item [A2.] There is a positive number $c$ such that $V^{-1}([0,c])$ is a compact subset of $U$.
 
\item [A3.] The set of  all critical points of $V$  in $V^{-1}([0,c])$ is equal to $V^{-1}(0)$. 
\end{itemize}
Adding the negative gradient of $V$ to (\ref{our:dyn:sys}), let us consider the following dynamical system on $U$:
\begin{equation}\label{our:new:dyn:sys}
\dot x = X(x) - \nabla V(x).
\end{equation}
Since $0$ is the minimum value of $V$, $\nabla V (x)= 0$ for all $x\in V^{-1}(0)$. Hence, the two vector fields $X$ and $X-\nabla V$ coincide on $V^{-1}(0)$.
\begin{theorem}[Theorem 2.1 in \cite{ChJiPe16}]\label{theorem:general}
Under assumptions A1 -- A3, every trajectory of (\ref{our:new:dyn:sys}) starting in $V^{-1}([0,c])$ remains in $V^{-1}([0,c])$ for all future time and asymptotically converges to the  set $V^{-1}(0)$ as $t \rightarrow \infty$. Furthermore, $V^{-1}(0)$ is an invariant set of   both  (\ref{our:dyn:sys}) and (\ref{our:new:dyn:sys}).
\end{theorem}
\begin{remark}
Theorem \ref{theorem:general}  still holds with 
the use of the following modified dynamics
\[
\dot x = X(x) - A(x) \nabla V(x)
\]
 instead of (\ref{our:new:dyn:sys}), where $A(x)$ is an $n\times n$ matrix-valued function with its symmetric part $(A(x) + A^T(x))$  positive definite at each $x\in\mathbb R^n$.  Refer to \cite{ChJiPe16} for more detail.
\end{remark}

We now  explain a strategy to design feedback integrators for nonholonomic mechanical systems. When  a nonholonomic mechanical system of the form \eqref{original:nonholo} is given, first  construct an extended system of the form \eqref{bam:qp} and further extend it to  Euclidean space if necessary. The final extended system corresponds to \eqref{our:dyn:sys}. 
Then, choose a Lyapunov function $V$ such that   assumptions A1 -- A3 in the above are satisfied and $V^{-1}(0)$ coincides with an invariant set of interest of the extended system. Usually, $V$ is chosen such that the invariant set $V^{-1}(0)$ is equal to or contained in  the intersection of a level set of the extended Hamiltonian $\tilde H$ and the nonholonomic constraint set $\mathcal C$. 
Adding the negative gradient of $V$ to the right side of \eqref{our:dyn:sys}, determines a feedback integrator of the form \eqref{our:new:dyn:sys}. Next, choose any initial point from $V^{-1}(0)$ and integrate the system  \eqref{our:new:dyn:sys} from this initial point using a general integration scheme such as Euler, Runge-Kutta or any other scheme.  Then, the  numerical trajectory will remain close to the set $V^{-1}(0)$. Rigorously speaking,  the  numerical trajectory converges to an attractor $\Lambda_h$ of the discrete-time dynamical system derived from the chosen one-step numerical integrator with uniform step size $h$, where $\Lambda_h$ converges to  $V^{-1}(0)$ as $h \rightarrow 0+$. Refer to \cite{ChJiPe16} for more details. The same strategy applies to symmetry-reduced nonholonomic mechanical systems.  In the following section, we apply this strategy to various  nonholonomic mechanical systems with or without symmetry to illustrate the design of feedback integrators for nonholonomic mechanical systems.

\section{Applications}

\subsection{The Suslov Problem on $\operatorname{SO}(3)$} 

\paragraph{Equations of Motion.}The equations of motion of a rigid body with a linear constraint are given by
\begin{subequations}\label{Suslov:original}
\begin{align}
\dot R &= R ({\mathbb I}^{-1}\Pi)^\wedge,\\
\dot \Pi &= \Pi \times {\mathbb I}^{-1}\Pi + \lambda \mathbb I e,\\
 \Pi \cdot e&=0,\label{Suslov:constraint}
\end{align}
\end{subequations}
where $(R,\Pi) \in \operatorname{SO}(3) \times \mathbb R^3$,  $e\in\mathbb R^3$ such that  
\begin{equation}\label{e:unit:Suslov}
e \cdot \mathbb Ie = 1,
\end{equation}
and $\mathbb I$ is the moment of inertia tensor of the rigid body. The hat map $\wedge: \mathbb R^3 \rightarrow \mathfrak {so}(3)$ is defined by
\[
(\Omega_1, \Omega_2, \Omega_3)^\wedge = \begin{bmatrix}
0 & -\Omega_3 & \Omega_2 \\
\Omega_3 & 0 & -\Omega_1 \\
-\Omega_2 & \Omega_1 & 0
\end{bmatrix}.
\]
 Equation \eqref{Suslov:constraint} is the nonholonomic constraint on the system; see p. 394 of \cite{Bl03}. The multiplier $\lambda$ is computed as  $\lambda = -e\cdot (\Pi \times \mathbb I^{-1}\Pi)$, and 
the  reduced Hamiltonian $h$ of \eqref{Suslov:original} is 
\begin{equation}\label{def:h:Suslov}
h( \Pi) =\frac{1}{2}\Pi\cdot \mathbb I^{-1}\Pi.
\end{equation}

\paragraph{Extension.} Let us now extend the system from the constraint set \eqref{Suslov:constraint} to the entire phase space $\operatorname{SO}(3) \times \mathbb R^3$.
The extended reduced Hamiltonian $\tilde h$ is given by
\begin{equation}\label{def:htilde:Suslov}
\tilde h (\Pi) = \frac{1}{2}\Pi \cdot \mathbb I^{-1}\Pi - \frac{1}{2} (\Pi \cdot e)^2,
\end{equation}
and the equations of motion of the corresponding extended system on $\operatorname{SO}(3) \times \mathbb R^3$ are given by
\begin{subequations}\label{Suslov:extended}
\begin{align}
\dot R &= R (\mathbb I^{-1}\Pi - (\Pi \cdot e) e)^\wedge,\\
\dot \Pi &= \Pi \times(\mathbb I^{-1}\Pi - (\Pi \cdot e) e)+ \tilde \lambda \mathbb I e,
\end{align}
\end{subequations}
where 
\begin{equation}\label{lambda:tilde:Suslov}
\tilde \lambda = -e\cdot (\Pi \times \mathbb I^{-1}\Pi).
\end{equation}
The system \eqref{Suslov:extended}  has at least two first integrals: the Hamiltonian $\tilde h$ and the constraint momentum map
\begin{equation}\label{def:J:Suslov}
J(\Pi) = \Pi \cdot e.
\end{equation}

\paragraph{Feedback Integrators.} Let us  implement a feedback integrator for the extended system \eqref{Suslov:extended} that preserves   the manifold $\operatorname{SO}(3)$ and the values of the Hamiltonian $\tilde h$ and the constraint momentum $J$. From here on, we regard the rotation matrix $R$ as a $3\times 3$ matrix in $\mathbb R^{3\times 3}$, thus extending the system \eqref{Suslov:extended}  further from $\operatorname{SO}(3) \times \mathbb R^3$  to $\mathbb R^{3\times 3} \times \mathbb R^3$.

Choose any number $\tilde h_0$ such that
\[
\tilde h_0>0.
\]
Let
\[
U = \{(R,\Pi) \in \mathbb R^{3\times 3} \times \mathbb R^3 \mid \det R >0\}.
\]
and define a Lyapunov function $V: U \rightarrow \mathbb R$ by
\begin{equation}\label{grad:V:Suslov}
 V(R,\Pi) = \frac{k_0}{4}\|R^TR - I\|^2  + \frac{k_1}{2} |J(\Pi )|^2 + \frac{k_2}{2}|\tilde h(\Pi ) -  \tilde h_0|^2,
\end{equation}
where $k_0$, $k_1$ and  $k_2$ are positive numbers, and $\|  \, \|$ is the standard trace norm on $\mathbb R^{3\times 3}$, i.e, $\| A\| = \sqrt{\operatorname{trace}(A^TA)}.$ It is easy to see that  
\begin{align*}
V^{-1}(0) 
&= \{ (R,\Pi) \in \mathbb R^{3\times 3} \times \mathbb R^3 \mid R \in \operatorname{SO}(3), J(\Pi) = 0,  h(\Pi) =   \tilde h_0 \}
\end{align*}
since $\tilde h (\Pi) = h(\Pi)$ when $J(\Pi)=0$.
The gradient vector of $V$ is computed as
\begin{align}\label{grad:tildeV:Suslov}
\nabla_R  V = k_0 R(R^TR-I), \quad \nabla_{\Pi}  V  = k_1  J(\Pi ) e + k_2 \Delta \tilde h ( \mathbb  I^{-1}\Pi - (\Pi \cdot e) e),
\end{align}
where
\begin{equation}\label{Deltas:Suslov}
\Delta \tilde h = \tilde h(\Pi) -  \tilde h_0
\end{equation}
\begin{lemma}\label{lemma1:Suslov}
\[
\langle (\nabla_RV,\nabla_\Pi V),  (R (\mathbb I^{-1}\Pi - (\Pi \cdot e) e)^\wedge,  \Pi \times(\mathbb I^{-1}\Pi - (\Pi \cdot e) e)+ \tilde \lambda \mathbb I e)\rangle = 0.
\]
\begin{proof}
It is straightforward from \eqref{e:unit:Suslov} and \eqref{lambda:tilde:Suslov}. Alternatively, we can use the argument that $V$ is a first integral of \eqref{Suslov:extended} since it consists of the first integrals $(R^TR-I)$, $J$ and $\tilde H$ of \eqref{Suslov:extended}, where it is easy to verify that $(R^TR-I)$ is a first integral of \eqref{Suslov:extended}.
\end{proof}
\end{lemma}

\begin{lemma}\label{lemma2:Suslov}
For any $c$ satisfying
\begin{equation}\label{inequality:c:Suslov}
0< c< \min\{{k_0}/{4}, k_2|\tilde h_0|^2/2\},
\end{equation}
 the set $V^{-1}([0,c])$ is compact. 
\begin{proof}
Take any $c$ satisfying  \eqref{inequality:c:Suslov}, and any  $(R,\Pi) \in V^{-1}([0,c])$. Then,
$\|R^TR -I\| \leq \sqrt{4c/k_0} <1$, so $\|R^TR\| <  \| I \| + 1 = 4$, which implies the set of $R$ such that $(R,\Pi) \in V^{-1}([0,c])$ is bounded. 
Since $(R,\Pi) \in V^{-1}([0,c])$, we have $|\Pi\cdot e| \leq \sqrt{2c/k_1}$ and $|\tilde h(\Pi) -  \tilde h_0| \leq\sqrt{2c/k_2}$. By the triangle inequality, $|\Pi\cdot \mathbb I^{-1}\Pi| \leq 2c/k_1 + 2|\tilde h_0| + 2\sqrt{2c/k_2}$. Thus, the set of $\Pi$ such that $(R,\Pi) \in V^{-1}([0,c])$ is bounded. Therefore, $V^{-1}([0,c])$ is compact, being bounded and closed.
\end{proof}
\end{lemma}

\begin{lemma}\label{lemma3:Suslov}

For any $c$ that satisfies \eqref{inequality:c:Suslov}, the set of critical point of $V$ in $V^{-1}([0,c])$ is equal to $V^{-1}(0)$. 
\begin{proof}
Take any $c$ that satisfies  \eqref{inequality:c:Suslov}.  Let $(R,\Pi)$ be any critical point of $V$ in $V^{-1}([0,c])$. Then, it satisfies
\begin{align}
0&=k_0R(R^TR-I), \label{crit:R:Suslov}\\
0&=k_1 J (\Pi)e + k_2 \Delta \tilde h ( \mathbb I^{-1} \Pi - (\Pi \cdot e) e). \label{crit:Pi:Suslov}
\end{align}
As shown in the proof of Lemma \ref{lemma2:Suslov}, $R$ satisfies $\|R^TR -I\| <1$, so it is invertible. Hence,  \eqref{crit:R:Suslov} implies $R\in \operatorname{SO}(3)$. Taking the inner product of \eqref{crit:Pi:Suslov} with $\mathbb Ie$ and using \eqref{e:unit:Suslov}, we obtain $J (\Pi )= 0$. Substitution of $J(\Pi) = \Pi \cdot e = 0$ in \eqref{crit:Pi:Suslov} yields
\[
\Delta \tilde h \mathbb I^{-1}\Pi = 0.
\]
Suppose $\Delta \tilde h \neq 0$. Then, the above equation implies $\mathbb I^{-1}\Pi=0$, thus $h(\Pi) = \frac{1}{2}\Pi\cdot \mathbb I^{-1}\Pi - \frac{1}{2}J(\Pi)^2 = 0$.  Then,
\[
V(R,\Pi) \geq \frac{k_2}{2}|\tilde h(\Pi) - \tilde h_0|^2 =  \frac{k_2}{2}| \tilde h_0|^2 >c,
\]
which contradicts $(R,\Pi) \in V^{-1}([0,c])$. Therefore, we must have $\Delta \tilde h = 0$. 

We have thus proved that the  critical point $(R,\Pi)$ is contained in $V^{-1}(0)$.  Since $0$ is the minimum value of $V$, every point in $V^{-1}(0)$ is a critical point of $V$. Therefore, the set of critical points of $V$ in $V^{-1}([0,c])$ equals $V^{-1}(0)$. 
\end{proof}
\end{lemma}
The feedback integrator system  for \eqref{Suslov:extended} with the function $V$ in \eqref{grad:V:Suslov} is given by
\begin{subequations}\label{Suslov:FI:tilde}
\begin{align}
\dot R &= R (\mathbb I^{-1}\Pi - (\Pi \cdot e) e)^\wedge - \nabla_RV,\\
\dot \Pi &= \Pi \times(\mathbb I^{-1}\Pi - (\Pi \cdot e) e)+ \tilde \lambda \mathbb I e - \nabla_\Pi V,
\end{align}
\end{subequations}
where $\nabla_{R}V$ and $\nabla_{\Pi} V$ are given in \eqref{grad:tildeV:Suslov}.

\begin{theorem}\label{theorem:FI:Suslov}
Let $c$ be any number that satisfies \eqref{inequality:c:Suslov}. Then,  every trajectory of \eqref{Suslov:FI:tilde} starting in $V^{-1}([0,c])$ remains in $V^{-1}([0,c])$ forward in time and  asymptotically converges to $V^{-1}(0)$ as $t\rightarrow \infty$. Moreover, $V^{-1}(0)$ is an invariant set of \eqref{Suslov:FI:tilde}.
\begin{proof}
The theorem follows from Lemmas \ref{lemma1:Suslov} -- \ref{lemma3:Suslov} in the above and Theorem 2.1 in \cite{ChJiPe16}.
\end{proof}
\end{theorem}

\begin{remark}
By Theorem \ref{theorem:reduced:first:integrals} and Corollary \ref{corollary:reduced:H:instead}, we can build another feedback integrator by using $h$ instead of $\tilde h$ in the construction of the Lyapunov function as follows:
\[
V (R,\Pi) = \frac{k_0}{4}\| R^TR-I\|^2 + \frac{k_1}{2}|J(\Pi)|^2 + \frac{k_2}{2}|h(\Pi) -  h_0|^2,
\]
where $h_0>0$.
The corresponding feedback integrator is in the same form as that in \eqref{Suslov:FI:tilde} but with the following  gradient vector of $V$:
\begin{align*}
\nabla_R  V = k_0 R(R^TR-I), \quad
 \nabla_{\Pi}  V  = k_1  J(\Pi ) e + k_2 ( h(\Pi) -  h_0)  \mathbb  I^{-1}\Pi.
\end{align*}
Theorem \ref{theorem:FI:Suslov} holds for this feedback integrator for any $c$ satisfying
\[
0< c< \min\{{k_0}/{4}, k_2|h_0|^2/2, k_1^2/2k_2\}.
\]
The proof is left to the reader.
\end{remark}
%
%
\begin{figure}[t]
\hspace{-15mm}
\includegraphics[scale = 0.4]{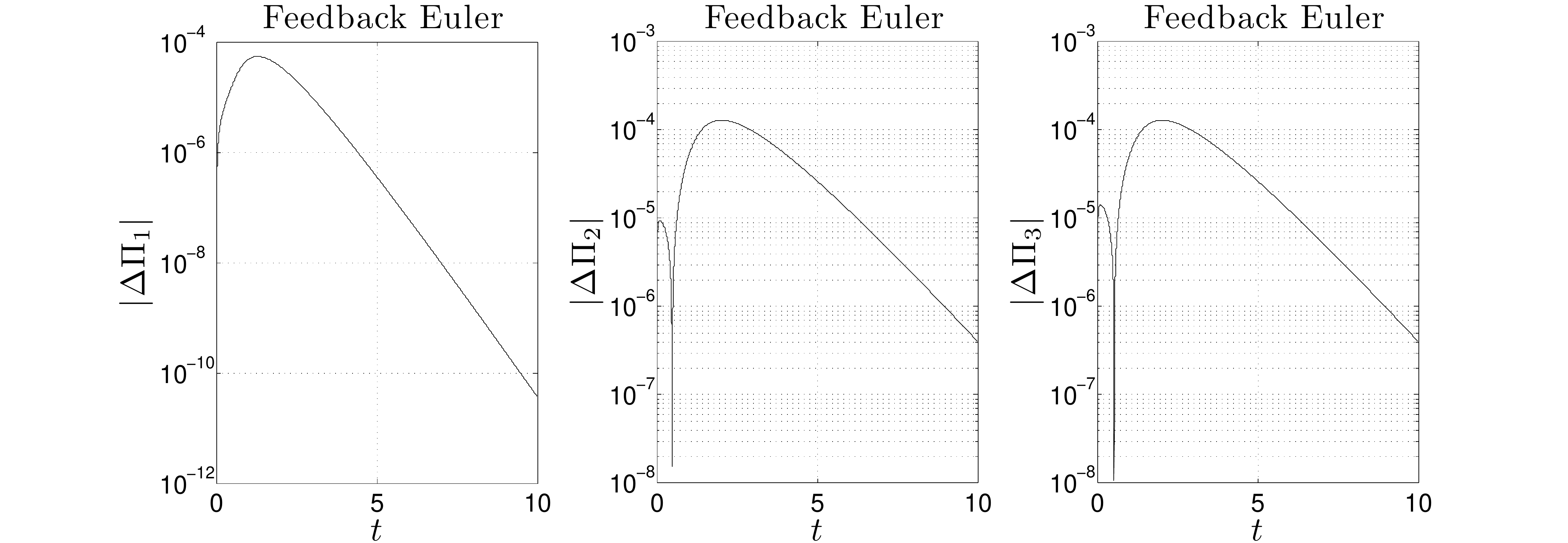}
\caption{The momentum error $\Delta \Pi (t) = \Pi(t) - \Pi_{\rm exact}(t)$, $0 \leq t\leq 10$, of the numerical solution $\Pi(t)$ of the Suslov system generated by  a feedback integrator with the Euler scheme with step size $\Delta t= 10^{-3}$ in comparison with the exact solution $ \Pi_{\rm exact}(t)$.}
\label{figure:Suslov_Pi_exact}
\end{figure}

\begin{figure}[t]
\hspace{-15mm}
\includegraphics[scale = 0.4]{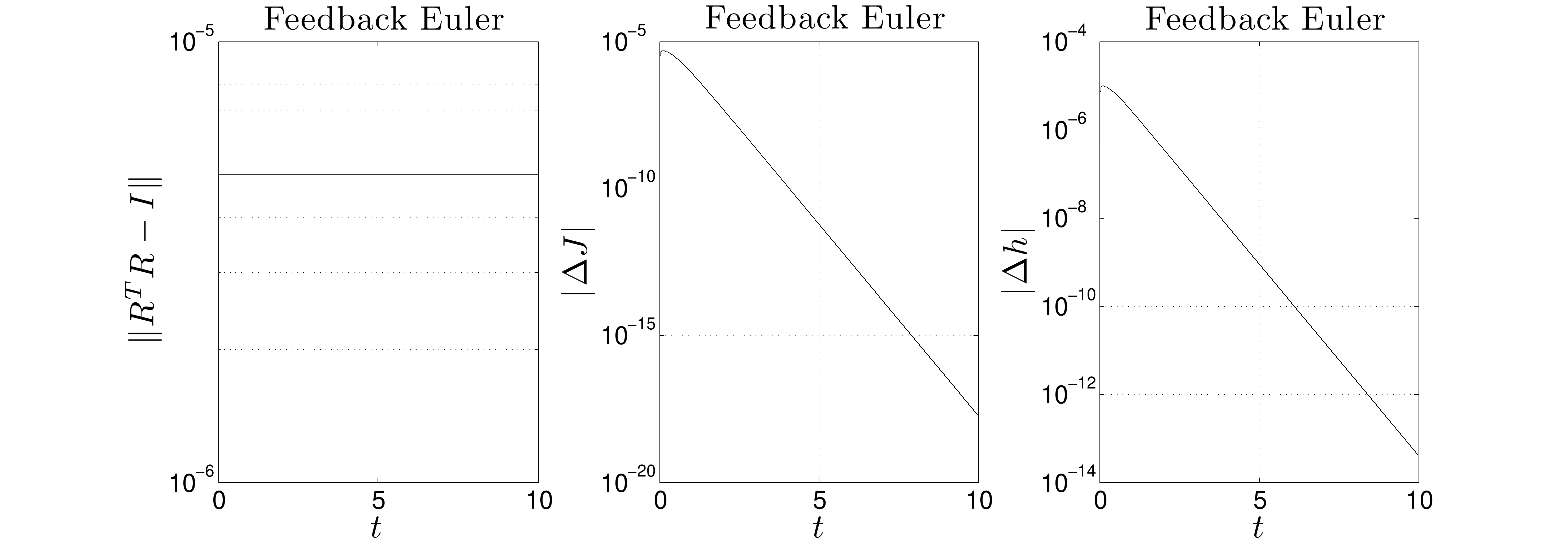}
\caption{The deviation $\|R(t)^TR(t) - I\|$ of $R(t)$ from the manifold $\operatorname{SO}(3)$, the nonholonomic constraint error $\Delta J (t)$, and the energy error $\Delta h (t)$ of the numerical solution $(R(t), \Pi(t))$, $0\leq t \leq 10$, of the Suslov system generated by  afeedback integrator with the Euler scheme with step size $\Delta t = 10^{-3}$.}
\label{figure:Suslov_integrals}
\end{figure}

\paragraph{Simulation.} Choose  parameter values as follows:
\[
\mathbb I =\begin{bmatrix}
1 & 0 & 0 \\
0 & 1 & 1 \\
0 & 1 & 2
\end{bmatrix}, \quad e = \frac{\mathbb I^{-1}a}{\sqrt{a\cdot \mathbb I^{-1}a}}
\]
where $a = (0,0,1)$, and the initial condition
\[
R(0) = I, \quad \Pi(0) = (0,1,1).
\]
The corresponding values of the momentum and  the Hamiltonian are 
\[
J(0) = 0, \quad \tilde h(0) = h(0) = 0.5.
\]
The exact solution $\Pi(t) = (\Pi_1(t), \Pi_2(t), \Pi_3(t))$  for the initial condition can be easily obtained as
\[
\Pi_1(t)= -\tanh t, \quad \Pi_2(t) = \Pi_3(t) = \frac{1}{\cosh t}.
\]
We take the integration time step size $\Delta t= 10^{-3}$ and  use the usual Euler method to integrate the feedback integrator system \eqref{Suslov:FI:tilde}  with the feedback gains 
\[
k_1=k_2=k_3 = 100
\]
over the time interval   $[0,10]$. Figure \ref{figure:Suslov_Pi_exact} shows the momentum error $\Delta \Pi = (\Delta \Pi_1,\Delta \Pi_2,\Delta \Pi_3)$ between the numerical solution and the exact solution. Figure \ref{figure:Suslov_integrals} shows the deviation $\| R^TR - I\|$ of the numerical solution $R(t)$ from the manifold $\operatorname{SO}(3)$, the constraint momentum error $|\Delta J|$ and the energy error $|\Delta h|$, where we use the original energy $h$ instead of the extended energy $\tilde h$ to show that  the conserved value of the original energy $h$ is well maintained. The simulation results demonstrate the excellent performance of the feedback integrator.

\subsection{The Knife Edge}
The system of a knife edge on an inclined plane  is an example in which the zero level set  $V^{-1}(0)$ of the Lyapunov function $V$ is not compact in spite of which we will see in a simulation that  the feedback integrator works well on the knife edge system. Hence, we expect with some caution that feedback integrators   perform well practically without the compactness assumption.

We follow the model of a knife edge on an inclined plane that appears in Section 1.6 of \cite{Bl03}. Let $\alpha >0$ denote the inclination angle of the plane and $(x,y)$ the position of the point of contact of the knife edge with respect to a fixed Cartesian coordinate system on the place (see Figure 1.6.1 in \cite{Bl03}). The angle $\varphi$ denotes the orientation angle of the knife edge with respect to the $xy$-plane. Let $m$ denote the mass of the knife and $J$ the moment of inertia of the knife edge about a vertical axis through its contact point. The gravitational acceleration constant is denoted by $g$.

The Hamiltonian $H$ of the knife edge system is given by
\[
H(q,p) = \frac{1}{2m}(p_x^2 + p_y^2) + \frac{1}{2J}p_\varphi^2 - mgx\sin\alpha,
\]
where $q= (x,y,\varphi)$ and $p = (p_x,p_y,p_\varphi)$. 
The equations of motion of the system  are given by
\begin{subequations}
\begin{align*}
\dot q &=\left (\frac{p_x}{m}, \frac{p_y}{m}, \frac{p_\varphi }{J} \right ),\\
\dot p &= (mg\sin\alpha, 0,0) + \lambda  m e,\\
0&=\frac{1}{\sqrt m}(p_x\sin \varphi-p_y\cos\varphi ),
\end{align*}
\end{subequations}
where
\begin{align*}
e=\frac{1}{\sqrt m} (\sin\varphi, -\cos\varphi,0).
\end{align*}

We now extend the system from the nonholonomic constraint set, $p_x\sin \varphi=p_y\cos\varphi$, to the entire phase space. The extended Hamiltonian $\tilde H$ is computed as
\begin{align*}
\tilde H(q,p) &=  \frac{1}{2m}(p_x^2 + p_y^2) + \frac{1}{2J}p_\varphi^2  - \frac{1}{2m}(p_x \sin\varphi - p_y\cos\varphi )^2 - mgx\sin \alpha  \\
&= \frac{1}{2m}(p_x\cos\varphi +p_y\sin\varphi)^2 + \frac{1}{2J}p_\varphi^2 - mgx\sin\alpha.
\end{align*} 
The equations of motion of the extended system are given by
\begin{subequations}\label{ext:sys:kinfe}
\begin{align}
\dot q &= \nabla_p \tilde H,\\
\dot p &= -\nabla_q\tilde H + \tilde \lambda m e,
\end{align}
\end{subequations}
where
\begin{align*}
\nabla_q\tilde H &= (-mg\sin\alpha,0,0),\\
\nabla_p\tilde H &= \left ( \frac{(p_x\cos\varphi + p_y\sin\varphi )\cos\varphi}{m},  \frac{(p_x\cos\varphi+p_y\sin\varphi )\sin\varphi}{m}, \frac{p_\varphi}{J} \right ),\\
\tilde \lambda &=\frac{1}{\sqrt{m}} \left ( -mg \sin\alpha\sin\varphi - \frac{p_\varphi(p_x\cos\varphi + p_y \sin\varphi)}{J} \right ). 
\end{align*}
The extended system \eqref{ext:sys:kinfe} has the following three first integrals: the Hamiltonian $\tilde H$ and  the two momentum maps $J_1(q,p)$ and $J_2(q,p)$ defined by
\[
J_1(q,p) = \frac{1}{\sqrt{m}}(p_x \sin\varphi - p_y \cos\varphi ), \quad J_2(q,p) = p_\varphi,
\]
the first of which comes from the nonholonomic constraint.

We now construct a feedback integrator for the extended system \eqref{ext:sys:kinfe}. Choose two  numbers $\tilde  H_0$ and $J_{20}$. Define a Lyapunov function $V$ by
\[
V(q,p) = \frac{k_1}{2}|J_1(q,p)|^2 + \frac{k_2}{2}|J_2(q,p) - J_{20}|^2 + \frac{k_3}{2}|\tilde H(q,p) - \tilde  H_0|^2,
\]
with $k_1,k_2,k_3>0$. Notice that the set
\[
V^{-1}(0) = \{  p_x \sin\varphi = p_y \cos\varphi , p_\varphi = J_{20}, H(q,p) = \tilde  H_0\}
\]
is not compact. 
The feedback integrator system corresponding to $V$ is given by
\begin{subequations}\label{FI:knife}
\begin{align}
\dot q &=\nabla_p\tilde H - \nabla_qV,\\
\dot p &= -\nabla_q \tilde H + \tilde \lambda m e- \nabla_pV,
\end{align}
\end{subequations}
where
\begin{align*}
\nabla_qV &= k_1 J_1(q,p) \nabla_qJ_1 + k_3 \Delta \tilde H\nabla_q\tilde H,\\
\nabla_pV&= k_1 J_1(q,p) \nabla_pJ_1 + k_2\Delta J_2 \nabla_pJ_2 + k_3\Delta \tilde H \nabla_p\tilde H
\end{align*}
and
\begin{align*}
\nabla_q J_1 &= \frac{1}{\sqrt{m}}(0,0,p_x\cos\varphi + p_y\sin\varphi), \,\,\, \nabla_pJ_1 = \frac{1}{\sqrt{m}}(\sin\varphi, -\cos\varphi,0), \\
 \nabla_pJ_2& = (0,0,1),\\
 \Delta \tilde H &= (\tilde H(q,p) - \tilde  H_0), \,\,\, \Delta J_2 = J_2(q,p)-J_{20}.
\end{align*}

\begin{theorem}\label{theorem:crit:knife}
The set of all critical points of $V$ equals $V^{-1}(0)$.
\begin{proof}
Choose an arbitrary critical point $(q,p)$ of $V$. Then, it satisfies
\begin{align}
0&= k_1 J_1(q,p) \nabla_qJ_1 + k_3 \Delta \tilde H\nabla_q\tilde H, \label{crit:A:knife}\\
0&= k_1 J_1(q,p) \nabla_pJ_1 + k_2\Delta J_2 \nabla_pJ_2 + k_3\Delta \tilde H \nabla_p\tilde H. \label{crit:B:knife}
\end{align}
Taking the inner product of \eqref{crit:A:knife} with $\nabla_q \tilde H$, we get $0=k_3 \Delta \tilde H (mg\sin  \alpha)^2$, and thus $\Delta \tilde H=0$ since $\sin\alpha >0$. Substituting $\Delta \tilde H=0$ in \eqref{crit:B:knife} and taking the inner product of \eqref{crit:B:knife} with $\nabla_pJ_2$, we get $\Delta J_2 = 0$.  Then,  \eqref{crit:B:knife}  reduces to $0= J_1(q,p) (\sin \varphi, -\cos\varphi,0)$. Hence, $J_1(q,p) = 0$. Thus, $(q,p) \in V^{-1}(0)$, implying that all critical points of $V$ are contained in $V^{-1}(0)$. Since $0$ is the minimum value of $V$, every point in $V^{-1}(0)$ is a critical point of $V$. Therefore, the set of all critical points of $V$ equals $V^{-1}(0)$.
\end{proof}
\end{theorem}
%
%
\begin{figure}[t]
\hspace{-15mm}
\includegraphics[scale = 0.37]{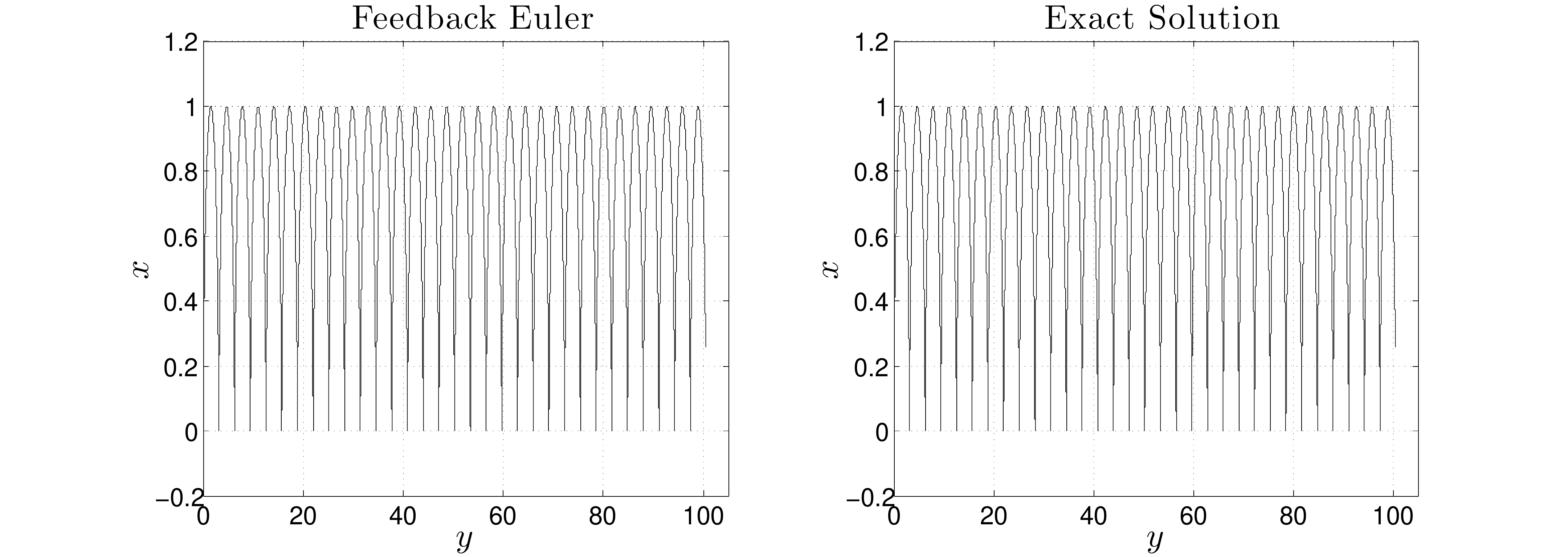}
\caption{The trajectory $(x(t), y(t))$, $0 \leq t \leq 200$, of the knife edge system generated by a feedback integrator with the Euler scheme with time step $\Delta t = 10^{-3}$, and the exact solution $(x(t), y(t)) = (1-\cos t,  ( t - \sin t )/2)$, $0 \leq t \leq 200$, where $y$ and $x$ are  swapped in the plots.}
\label{figure:knife_cycloid}
\end{figure}

\begin{figure}[t]
\hspace{-15mm}
\includegraphics[scale = 0.37]{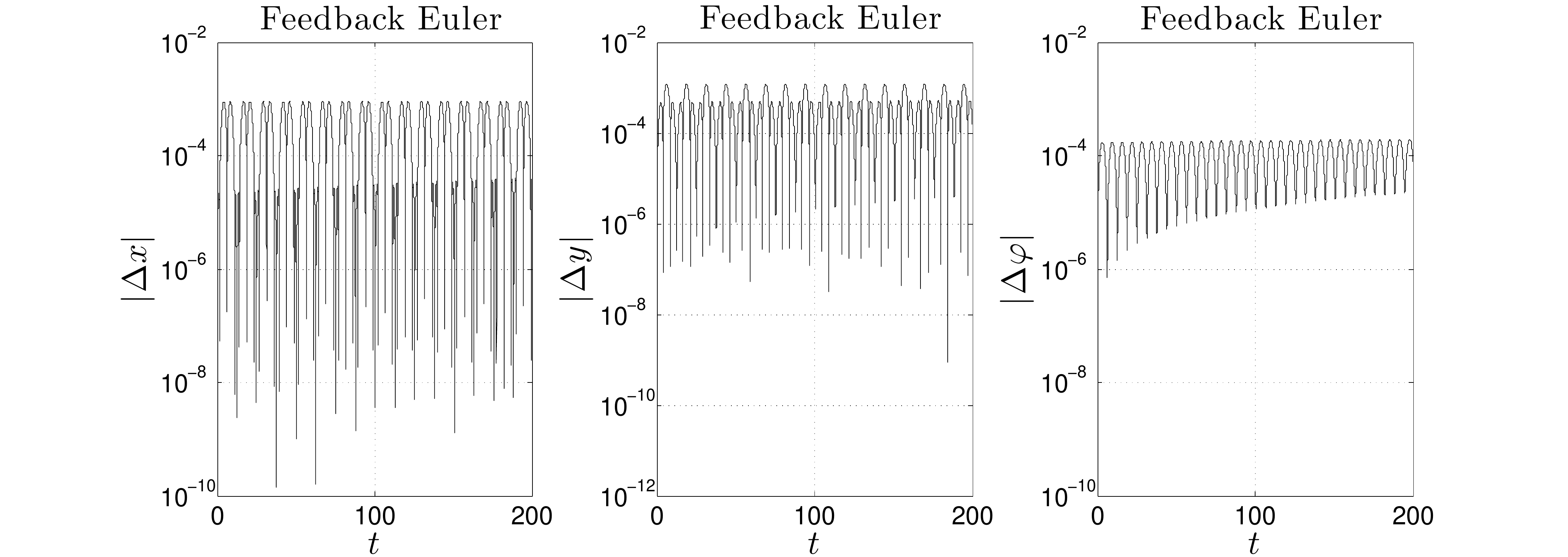}
\caption{The error of the trajectory $(x(t), y(t), \varphi (t))$, $0 \leq t \leq 200$,  of the knife edge system generated  by a feedback integrator with the Euler scheme with time step $\Delta t = 10^{-3}$ in comparison with the exact solution $(1-\cos t,  ( t - \sin t )/2, t)$, $0 \leq t \leq 200$.}
\label{figure:knife_exact}
\end{figure}

\begin{figure}[t]
\hspace{-15mm}
\includegraphics[scale = 0.36]{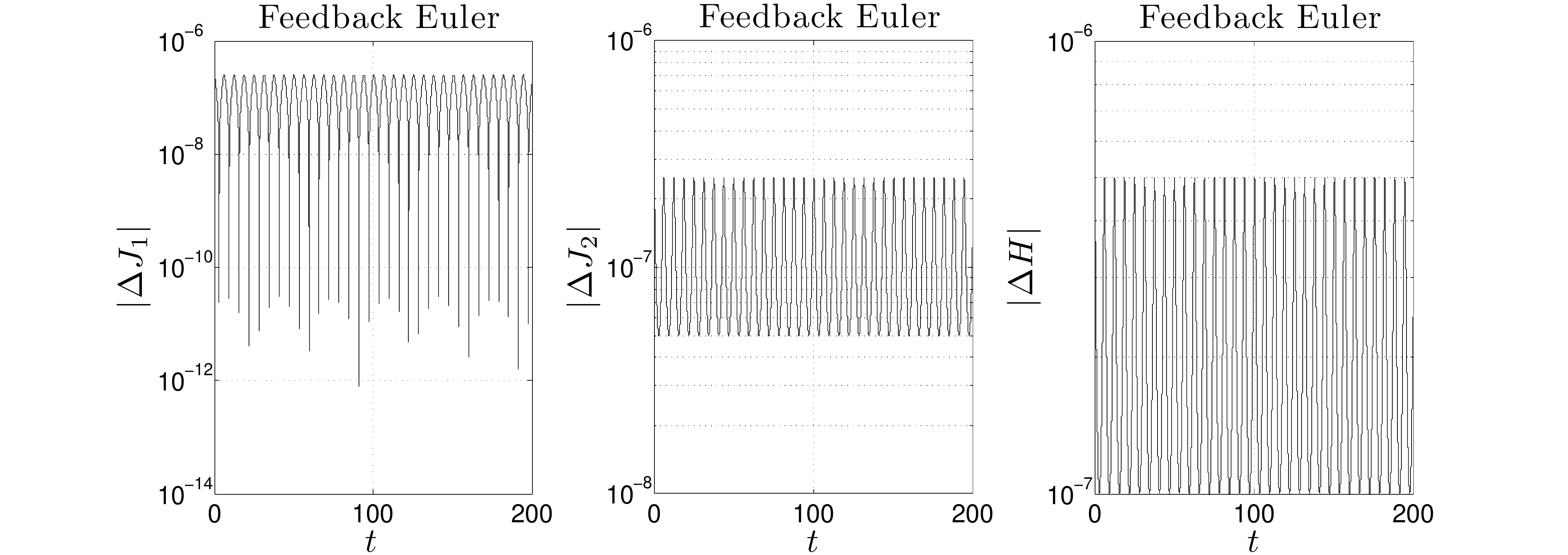}
\caption{The nonholonomic momentum error $\Delta J_1$, the momentum error $\Delta J_2 = \Delta p_\varphi$, and the energy error $\Delta H$ of the numerical solution  of the knife edge system generated by  a feedback integrator with the Euler scheme with step size $\Delta t = 10^{-3}$ over the time interval $[0,200]$.}
\label{figure:knife_integrals}
\end{figure}
\begin{remark}
We can design another feedback integrator for the extended system  \eqref{ext:sys:kinfe}  by using $H$ instead of $\tilde H$ in the construction of the Lyapunov function as follows:
\[
V(q,p) = \frac{k_1}{2}|J_1(q,p)|^2 + \frac{k_2}{2}|J_2(q,p) - J_{20}|^2 + \frac{k_3}{2}| H(q,p) -   H_0|^2.
\]
The corresponding feedback integrator is in the same form as that in \eqref{FI:knife} but with the following  gradient vector $\nabla V = (\nabla_q V, \nabla_p V)$ of $V$:
\begin{align*}
\nabla_q V &= k_1 J_1(q,p) \nabla_qJ_1 + k_3 \Delta  H(-mg\sin\alpha,0,0), \label{crit:A:knife}\\
\nabla_pV &= k_1 J_1(q,p) \nabla_pJ_1 + k_2\Delta J_2 \nabla_pJ_2 + k_3\Delta H p, 
\end{align*}
where $\Delta H = H(q,p) - H_0$.  Theorem \ref{theorem:crit:knife} also holds for this new Lyapunov function $V$, whose proof is left to the reader.
\end{remark}
\paragraph{Simulation.} Choose the parameter values 
\[
m = J = g = 1, \quad \alpha = \frac{\pi}{6}.
\]
and the initial conditions
\[
x(0) = y(0) = \varphi (0) = p_x (0) = p_y (0) = 0, \quad p_{\varphi} (0) = 0.5.
\]
The corresponding values of the momentum and  the Hamiltonian are 
\[
J_1(0) = 0, \quad J_2(0) = 0.5, \quad \tilde H(0) = H(0) = 0.125.
\]
The exact solution $(x(t), y(t))$  for these initial conditions are given by
\[
x(t) =1-\cos t, \quad  y(t) = \frac{1}{2} ( t - \sin t ), \quad \varphi(t) =t,
\]
so that $(x(t),y(t))$ undergoes a cycloid motion; refer to Section 1.6 of \cite{Bl03} for the derivation of the exact solution.

We take the integration time step size $\Delta t = 10^{-3}$ and  use the usual Euler method to integrate the feedback integrator system \eqref{theorem:crit:knife}  with the feedback gains 
\[
k_1=k_2=k_3 = 1000
\]
over the time interval   $[0,200]$. Figure \ref{figure:knife_cycloid} compares the motion undergone by the numerical solution $(x(t), y(t))$ with the cycloid motion of the exact solution. Figure \ref{figure:knife_exact} provides the plots of the errors $|\Delta x|$, $|\Delta y|$, and $|\Delta \varphi |$ between the numerical  solution and the exact solution. Figure \ref{figure:knife_integrals} shows the plots of the  momenta errors $|\Delta J_1|$ and $|\Delta J_2|$ and the energy error $|\Delta h|$, where  we use the original energy $h$ instead of the extended energy $\tilde h$ to show that  the conserved value of the original energy $h$ is well maintained. 
The simulation results demonstrate the efficacy of the feedback integrator even in the absence of compactness of the set $V^{-1}(0)$.

\subsection{The Chaplygin Sleigh}
Consider the Chaplygin sleigh system in Figure 1.7.1 in \cite{Bl03}. 
The configuration space is the Lie group $ \operatorname{SE}(2)$; refer to section 14.6 of \cite{MaRa02} for the Lie group $\operatorname{SE}(2)$ and its dual Lie algebra $\mathfrak{se}(2)^* = \mathbb R\times \mathbb R^2$. We use  $q = (\theta, x,y)$ for coordinates on  $\operatorname{SE}(2)$ and $p = (p_\theta, p_x,p_y)$ for the corresponding momentum.

The Hamiltonian $H$ is
\[
H(q,p) =\frac{1}{2I} (p_\theta - a(-p_x\sin\theta + p_y \cos\theta))^2 +  \frac{1}{2m}(p_x^2 + p_y^2),
\]
where $m$ is the mass,  $I$  the moment inertia of the sleigh, and $a$ is the distance from the contact point of the knife edge to the center of mass of the sleigh.
The nonholonomic constraint is
\[
 \mathcal C = \{ ma p_\theta - (I+ma^2) (-p_x \sin\theta + p_y \cos\theta) = 0 \}.
\]
Since both $H$ and $\mathcal C$ are $\operatorname{SE}(2)$-invariant under left action, they can be reduced to the dual Lie algebra $\mathfrak {se}(2)^* = \mathbb R \times \mathbb R^2 = \mathbb R^3$. According to equation (14.6.16) in \cite{MaRa02}, the minus Lie-Poisson bracket on $\mathfrak{se}(2)^*$ is 
\begin{equation}\label{Poisson:bracket:se2}
\{ F, G\}_{\mathfrak{se}(2)^*}(\mu,\alpha) = - \left ( \frac{\partial F}{\partial \mu} \mathbb J\alpha \cdot \nabla_\alpha G - \frac{\partial G}{\partial \mu} \mathbb J\alpha \cdot \nabla_\alpha F\right )
\end{equation}
for  any $(\mu, \alpha) \in \mathbb R \times \mathbb R^2 \simeq \mathfrak{se}(2)^*$ and  any $F,G \in C^\infty( \mathfrak{se}(2)^*)$, where
\[
\mathbb J = \begin{bmatrix}
0 & 1 \\ -1 & 0
\end{bmatrix}.
\]
 Let $\Pi = (\Pi_\theta, \Pi_x, \Pi_y)$ be the coordinates for $\mathfrak {se}(2)^* = \mathbb R^3$ such that
\[
 (\Pi_\theta, \Pi_x, \Pi_y) = (p_\theta, p_x \cos\theta +  p_y \sin\theta , -p_x\sin \theta + p_y \cos\theta).
\]
The reduced Hamiltonian $h$ and the reduced constraint set ${\mathcal C}_{\rm red}$ are given by
\[
h(\Pi) = \frac{1}{2I}  (\Pi_\theta - a \Pi_y)^2 + \frac{1}{2m} (\Pi_x^2 + \Pi_y^2)      
\]
and
\[
{\mathcal C}_{\rm red} = \{ma \Pi_\theta - (I+ma^2) \Pi_y = 0\} = \{ \langle \Pi, e\rangle = 0\} ,
\]
where
\[
e = \frac{1}{\sqrt{mI (I + ma^2)}} (ma, 0, -(I + ma^2)).
\]
Notice that the vector $e$ has unit length with respect to the locked inertia tensor $\mathbb I$  that is given by
\[
\mathbb I = \begin{bmatrix}
I + ma^2 & 0 & ma \\
0 & m & 0 \\
ma & 0 & m
\end{bmatrix}.
\]
Hence,
\begin{equation}\label{e:flat:Chap}
e^\flat = \frac{-mI}{\sqrt{mI (I + ma^2)}} (0, 0, 1).
\end{equation}
The extended reduced Hamiltonian $\tilde h$ is computed as
\begin{equation}\label{h:tilde:Chap}
\tilde h (\Pi) = \frac{1}{2(I+ma^2)}\Pi_\theta^2 + \frac{1}{2m}\Pi_x^2.
\end{equation}
Compute
\begin{equation}\label{lambda:Chap}
\tilde \lambda = \{ \tilde h, \langle p,e\rangle \}_- = \frac{-1}{\sqrt{mI (I + ma^2)}} (\Pi_\theta - a \Pi_y)\Pi_x,
\end{equation}
where the bracket formula can be found in equation (14.6.16) in \cite{MaRa02}. The reduced extended system, corresponding, in the general theory section to equations \eqref{new:red:bam} is  given in general form by
\begin{subequations}\label{general:form:SE2:eqn}
\begin{align}
(\dot \theta, \dot x, \dot y) &=  \left (\frac{\partial \tilde h}{\partial \Pi_\theta}, \frac{\partial \tilde h}{\partial \Pi_x} \cos \theta - \frac{\partial \tilde h}{\partial \Pi_y} \sin \theta, \frac{\partial \tilde h}{\partial \Pi_x} \sin \theta + \frac{\partial \tilde h}{\partial \Pi_y} \cos \theta \right ),\\
(\dot \Pi_\theta, \dot \Pi_x, \dot \Pi_y) &= \left ( -\Pi_y \frac{\partial \tilde h}{\partial \Pi_x} + \Pi_x \frac{\partial \tilde h}{\partial \Pi_y}, \Pi_y \frac{\partial \tilde h}{\partial \Pi_\theta}, -\Pi_x \frac{\partial \tilde h}{\partial \Pi_\theta} \right ) + \tilde \lambda e^\flat,  
\end{align}
\end{subequations}
which reduces, with substitution of  \eqref{e:flat:Chap},  \eqref{h:tilde:Chap} and \eqref{lambda:Chap}, to  
\begin{subequations}\label{extended:system:Chap}
\begin{align}
(\dot \theta, \dot x, \dot y) &= \left ( \frac{\Pi_\theta}{I+ma^2}, \frac{\Pi_x \cos\theta}{m},  \frac{\Pi_x \sin\theta}{m} \right ), \label{extended:system:q:Chap}\\
(\dot \Pi_\theta, \dot \Pi_x, \dot \Pi_y) &=  \left (-\frac{\Pi_x\Pi_y}{m}, \frac{\Pi_\theta \Pi_y}{I+ma^2}, -\frac{a\Pi_x \Pi_y}{I+ma^2} \right ). \label{extended:system:Pi:Chap}
\end{align}
\end{subequations}
We now construct a feedback integrator for the extended system \eqref{extended:system:Chap}. Choose any number $\tilde h_0$ such that
\[
\tilde h_0 >0. 
\]
Define a Lyapunov function $V$ on $\mathbb R^3$ by
\[
V(\Pi) = \frac{k_1}{2}|J(\Pi)|^2 + \frac{k_2}{2}|\tilde h( \Pi) - \tilde h_0|^2
\]
where $k_1, k_2>0$ and
\[
J(\Pi) = \langle \Pi, e\rangle.
\]
The function $V$  satisfies
\[
V^{-1}(0) =\{ \Pi \in \mathbb R^3 \mid \Pi \in \mathbb R^3 \mid J(\Pi) = 0,  h(\Pi) = \tilde h_0\}
\]
since $\tilde h(\Pi) = h(\Pi)$ when $J(\Pi) = 0$. The feedback integrator corresponding to the function $V$ is given by
\begin{subequations}\label{FI:Chap}
\begin{align}
(\dot \theta, \dot x, \dot y) &= \left ( \frac{\Pi_\theta}{I+ma^2}, \frac{\Pi_x \cos\theta}{m},  \frac{\Pi_x \sin\theta}{m} \right ), \label{FI:q:Chap}\\
(\dot \Pi_\theta, \dot \Pi_x, \dot \Pi_y) &=  \left (-\frac{\Pi_x\Pi_y}{m}, \frac{\Pi_\theta \Pi_y}{I+ma^2}, -\frac{a\Pi_x \Pi_y}{I+ma^2} \right ) - \nabla_{\Pi}V, \label{FI:Pi:Chap}
\end{align}
\end{subequations}
where
\begin{equation}\label{nabla:Pi:V:Chap}
\nabla_{\Pi}V  =k_1 J(\Pi) e + k_2 (\tilde h(\Pi) - \tilde h_0) \left ( \frac{\Pi_\theta}{I + ma^2}, \frac{\Pi_x}{m},0\right).
\end{equation}
Notice that the subsystem \eqref{FI:Pi:Chap} is the essential part of the feedback integrator that is not affected by the other part \eqref{FI:q:Chap}. 
\begin{theorem}\label{theorem:FI:Chap}
For any number $c$ satisfying $0< c< k_2|\tilde h_0|^2/2$,  every trajectory of \eqref{FI:Chap} starting in $\operatorname{SE}(2)\times V^{-1}([0,c])$ remains in $\operatorname{SE}(2)\times V^{-1}([0,c])$ for all future time and asymptotically converges to $\operatorname{SE}(2)\times V^{-1}(0)$ as $t \rightarrow\infty$. Moreover, $\operatorname{SE}(2)\times V^{-1}(0)$ is an invariant set of \eqref{FI:Chap}.
\begin{proof}
Since both $J$ and $\tilde h$ are first integrals of \eqref{extended:system:Chap}, the function $V$ is also a first integral of \eqref{extended:system:Chap}.
Take any number $c$ such that $0 < c<k_2|\tilde h_0|^2/2$. It is easy to show that  the set $V^{1}([0,c])$ is a compact subset of $\mathbb R^3$. Take any critical point $\Pi$ of $V$ in $V^{-1}([0,c])$. Then, it satisfies
\[
0=k_1 J(\Pi) e + k_2 (\tilde h(\Pi) - \tilde h_0) \left ( \frac{\Pi_\theta}{I + ma^2}, \frac{\Pi_x}{m},0\right).
\]
From the third component of the above vector equation, it follows that $J_1(\Pi) = 0$. If $ \tilde h(\Pi) \neq \tilde h_0$, then $\Pi_\theta = \Pi_x = 0$, thus implying that $\tilde h(\Pi) = 0$, which would imply $V(\Pi) \geq k_2 |\tilde h_0|^2/2>c$, contradicting $\Pi \in V^{-1}([0,c])$. Hence, $\tilde h(\Pi) = \tilde h_0$. Thus, $\Pi \in V^{-1}(0)$, implying that every critical point of $V$ in $V^{-1}([0,c])$ belongs to $V^{-1}(0)$. Since $0$ is the minimum value of $V$, every point of $V^{-1}(0)$ is a critical point of $V$. Hence, by Theorem 2.1 in \cite{ChJiPe16}, every trajectory of \eqref{FI:Pi:Chap} starting in $V^{-1}([0,c])$ remains in $ V^{-1}([0,c])$ for all future time and asymptotically converges to $V^{-1}(0)$ as $t \rightarrow\infty$. Also, $V^{-1}(0)$ is an invariant set of \eqref{FI:Pi:Chap}.  Therefore, the theorem holds.
\end{proof}
\end{theorem}

We now design another feedback integrator for the extended system \eqref{extended:system:Chap} by embedding the factor $\operatorname{S}^1$ of $\operatorname{SE}(2)$ into $\mathbb R^{2\times 2}$ via the isomorphism between $\operatorname{S}^1$  and $\operatorname{SO}(2) = \{ R \in \mathbb R^{2\times 2} \mid R^T R =I\}$ given by
\[
\theta \mapsto \begin{bmatrix}
\cos\theta & - \sin\theta\\
\sin\theta & \cos \theta
\end{bmatrix}.
\]  Via the isomorphism, the general form of equations \eqref{general:form:SE2:eqn} can be written as
\begin{subequations}\label{general:form:SE2:w:SO2}
\begin{align}
\dot R &= -R {\mathbb J} \frac{\partial \tilde h}{\partial \Pi_\theta},\quad
\begin{bmatrix}
\dot x \\ \dot y
\end{bmatrix} = R \begin{bmatrix}
\displaystyle \frac{\partial \tilde h}{\partial \Pi_x} \\ \displaystyle  \frac{\partial \tilde h}{\partial \Pi_y}
\end{bmatrix},\\
(\dot \Pi_\theta, \dot \Pi_x, \dot \Pi_y) &= \left ( -\Pi_y \frac{\partial \tilde h}{\partial \Pi_x} + \Pi_x \frac{\partial \tilde h}{\partial \Pi_y}, \Pi_y \frac{\partial \tilde h}{\partial \Pi_\theta}, -\Pi_x \frac{\partial \tilde h}{\partial \Pi_\theta} \right ) + \tilde \lambda e^\flat,  
\end{align}
\end{subequations}
where $R \in \operatorname{SO}(2)$ and 
\[
{\mathbb J} = \begin{bmatrix}
0 & 1 \\-1 & 0
\end{bmatrix}.
\]
With substitution of \eqref{e:flat:Chap},  \eqref{h:tilde:Chap} and \eqref{lambda:Chap}, the general form of system \eqref{general:form:SE2:w:SO2} becomes 
\begin{subequations}\label{extended:SE2:SO2}
\begin{align}
\dot R &= -\frac{\Pi_\theta}{I+ma^2} R {\mathbb J},\\
\begin{bmatrix}
\dot x \\ \dot y
\end{bmatrix} &= R \begin{bmatrix}
\displaystyle  \frac{\Pi_x}{m} \\ 0
\end{bmatrix},\\
(\dot \Pi_\theta, \dot \Pi_x, \dot \Pi_y) &=  \left (-\frac{\Pi_x\Pi_y}{m}, \frac{\Pi_\theta \Pi_y}{I+ma^2}, -\frac{a\Pi_x \Pi_y}{I+ma^2} \right ),  
\end{align}
\end{subequations}
which is equivalent to \eqref{general:form:SE2:eqn}. We now treat the matrix $R$ as a $2\times 2$ matrix, extending the system \eqref{extended:SE2:SO2} further to $\mathbb R^{2 \times 2}\times \mathbb R^2 \times \mathbb R^3$. Choose any number $\tilde h_0>0$ and consider the Lyapunov function
\[
V(R,\Pi) = \frac{k_0}{2}\| R^TR - I\|^2+\frac{k_1}{2}|J_1(\Pi)|^2 + \frac{k_2}{2}|\tilde h( \Pi) - \tilde h_0|^2.
\]
The corresponding feedback integrator is computed as
\begin{subequations}\label{FI:extended:SE2:SO2}
\begin{align}
\dot R &= -\frac{\Pi_\theta}{I+ma^2} R {\mathbb J} - \nabla_RV, \label{FI:extended:R:SE2:SO2}\\
\begin{bmatrix}
\dot x \\ \dot y
\end{bmatrix} &= R \begin{bmatrix}
\displaystyle  \frac{\Pi_x}{m} \\ 0
\end{bmatrix}, \label{FI:extended:xy:SE2:SO2}\\
(\dot \Pi_\theta, \dot \Pi_x, \dot \Pi_y) &=  \left (-\frac{\Pi_x\Pi_y}{m}, \frac{\Pi_\theta \Pi_y}{I+ma^2}, -\frac{a\Pi_x \Pi_y}{I+ma^2} \right ) - \nabla_{\Pi}V,  
\end{align}
\end{subequations}
where $\nabla_R V$ is given by
\[
\nabla_R V = -k_0 R(R^TR -I)
\]
 and $\nabla_\Pi V$ is given in \eqref{nabla:Pi:V:Chap}. It is easy to prove that a theorem similar to Theorem \ref{theorem:FI:Chap} holds of the new feedback integrator \eqref{FI:extended:SE2:SO2}, whose proof is left to the reader. 
 
%

\subsection{The Vertical Rolling Disk}
We follow the model of the vertical rolling disk described in Section 1.4 of \cite{Bl03}. See Figure 1.4.1 therein with the replacement of   $\theta$  and $\varphi$ with  $\psi$ and $\theta$, respectively.  The configuration space is the Lie group $G = \operatorname{SE}(2) \times \operatorname{S}^1$. As for coordinates,  $(\theta, x,y)$ is used for $\operatorname{SE}(2)$ and $\psi$ for $ \operatorname{S}^1$.  The coordinates $(p_\theta, p_x,p_y,p_\psi )$ are for the corresponding conjugate momenta. 
The  Hamiltonian $H$ of the system is
\[
H(q,p) = \frac{1}{2J}p_\theta^2 +  \frac{1}{2m}(p_x^2 + p_y^2) + \frac{1}{2I}p_\psi^2 ,
\]
where $q = ( \theta,x,y,\psi)$ and $p = (p_\theta, p_x,p_y,p_\psi )$. The parameter $J$ is the moment of inertial about an axis in the plane of the disk,  $m$ is the mass of the disk, and $I$ is the moment of the inertia of the disc about the axis perpendicular to the plane of the disk. The set of nonholonomic constraints is given by
\begin{align*}
\mathcal C &= \left \{ 
\frac{1 }{m}p_x- \frac{R}{I} p_\psi \cos\theta= 0, \quad \frac{1 }{m}p_y- \frac{R}{I} p_\psi \sin\theta= 0 \right \} \\
&= \{ -p_x\sin \theta + p_y \cos\theta = 0,  I (p_x \cos\theta +  p_y \sin\theta ) - mRp_\psi = 0 \}.
\end{align*}
Since both $H$ and $\mathcal C$ are $G$-invariant under left action, they can be reduced to the dual Lie algebra ${\mathfrak g}^* = \mathfrak{se}(2)^* \times \mathbb  R=\mathbb R^4$. Let $\Pi = (\Pi_\theta, \Pi_x, \Pi_y, \Pi_\psi)$ be the coordinates for ${\mathfrak g}^* = \mathbb R^4$ such that
\[
 (\Pi_\theta, \Pi_x, \Pi_y, \Pi_\psi) = (p_\theta, p_x \cos\theta +  p_y \sin\theta , -p_x\sin \theta + p_y \cos\theta, p_\psi).
\]
The reduced Hamiltonian $h$  and the reduced set of nonholonomic constraints ${\mathcal C}_{\rm red}$ are  given by
\[
h(\Pi) = \frac{1}{2J}\Pi_\theta^2 + \frac{1}{2m}(\Pi_x^2 + \Pi_y^2) + \frac{1}{2I}\Pi_\psi^2
\]
and
\begin{align*}
{\mathcal C}_{\rm red} &= \{ \Pi_y = 0, I\Pi_x - mR\Pi_\psi=0 \}\\
&= \{ \langle \Pi, e_1\rangle = 0,  \langle \Pi, e_2\rangle = 0\},
\end{align*}
where
\[
e_1 = \frac{1}{\sqrt m}(0,0 , 1, 0), \quad e_2 =\frac{1}{\sqrt{mI(I+mR^2)}}(0, I  , 0, -mR).
\]
Notice that the two vectors $e_1$ and $e_2$ are orthonormal with respect to the locked inertia tensor $\mathbb I = \operatorname{diag}(J, m,m,I)$. Then, $e_1^b$ and $e_2^b$ are computed as
\[
e_1^\flat = (0,0 , \sqrt m, 0), \quad e_2^\flat =\frac{1}{\sqrt{mI(I+mR^2)}}(0, mI  , 0, -mRI).
\]
The extended Hamiltonian $\tilde h$ is given by
\[
\tilde h (\Pi)  = \frac{1}{2J}\Pi_\theta^2 + \frac{(R\Pi_x + \Pi_\psi)^2}{2(I+mR^2)}.
\]
The equations of motion of the corresponding system are 
\begin{subequations}\label{extended:eqn:vertical:disk}
\begin{align}
(\dot\theta, \dot x, \dot y, \dot \psi) &= \left (\frac{\Pi_\theta}{J}, \frac{R(R\Pi_x + \Pi_\psi)\cos\theta}{I+mR^2}, \frac{R(R\Pi_x + \Pi_\psi)\sin\theta}{I+mR^2}, \frac{R\Pi_x + \Pi_\psi}{I+mR^2} \right ),\\
(\dot \Pi_\theta, \dot \Pi_x,\dot \Pi_y,\dot \Pi_\psi)&=\left  (-\frac{R(R\Pi_x + \Pi_\psi)\Pi_y}{I+mR^2}, \frac{\Pi_\theta \Pi_y}{J}, -\frac{\Pi_\theta \Pi_x}{J},0 \right ) + \sum_{i=1}^2 \tilde \lambda_i e_i^\flat ,
\end{align}
\end{subequations}
where 
\[
\tilde \lambda_1 = \frac{1}{J\sqrt m}\Pi_\theta \Pi_x, \quad \tilde \lambda_2 = -\frac{I}{J\sqrt{mI(I+mR^2)}}\Pi_\theta\Pi_y.
\]
The extended system has the following three first integrals:  $\tilde h$, $\Pi_y$ and $I\Pi_x - mR\Pi_\psi$.  
We now further modify the extended system \eqref{extended:eqn:vertical:disk} outside $\mathcal C_{\rm red}$ while maintaining the constancy of motion of    $\tilde h$, $\Pi_y$ and $I\Pi_x - mR\Pi_\psi$, and maintaining the values of  $\tilde h$, $\Pi_y$ and $I\Pi_x - mR\Pi_\psi$ on $\mathcal C_{\rm red}$. For such a modification,  we substitute $\Pi_y = 0$ and $I\Pi_x - mR\Pi_\psi=0$ in \eqref{extended:eqn:vertical:disk}  to obtain
\begin{subequations}\label{further:extended:eqn:vertical:disk}
\begin{align}
(\dot\theta, \dot x, \dot y, \dot \psi) &= \left (\frac{\Pi_\theta}{J}, \frac{\Pi_x\cos\theta}{m}, \frac{\Pi_x\sin\theta}{m}, \frac{\Pi_\psi}{I} \right ),\\
(\dot \Pi_\theta, \dot \Pi_x,\dot \Pi_y,\dot \Pi_\psi)&=(0,0,0,0). 
\end{align}
\end{subequations}
It is easy to verify that the three functions  $\tilde h$, $\Pi_y$ and $I\Pi_x - mR\Pi_\psi$ are first integrals of \eqref{further:extended:eqn:vertical:disk} on the entire phase space as required. Moreover, all four momentum variables $\Pi_\theta$, $ \Pi_x$, $\Pi_y$, and $\Pi_\psi$ are first integrals of \eqref{further:extended:eqn:vertical:disk}, too.  Although it is trivial to integrate \eqref{further:extended:eqn:vertical:disk}, let us build a feedback integrator for \eqref{further:extended:eqn:vertical:disk} with the four conserved momentum variables.  Choose any numbers $\Pi_{\theta 0}$, $\Pi_{x0}$, $\Pi_{y0}$ and $\Pi_{\psi 0}$ such that 
$(\Pi_{\theta 0}, \Pi_{x0},\Pi_{y0},\Pi_{\psi 0}) \in {\mathcal C}_{\rm red}$, i.e.
\[
\Pi_{y0} = 0, \quad I\Pi_{x0} - R\Pi_{\psi 0} =0.
\] Consider the Lyapunov function
\[
V(\Pi)  = \frac{k_1}{2}|\Pi_\theta - \Pi_{\theta 0}|^2 + \frac{k_2}{2}|\Pi_x - \Pi_{x0}|^2 + \frac{k_3}{2}|\Pi_y|^2 + \frac{k_4}{2}| \Pi_\psi - \Pi_{\psi 0}|^2,
\]
where $k_i >0$ for all $i = 1, \ldots,4$. Then, the corresponding feedback integrator for \eqref{further:extended:eqn:vertical:disk}  is given by
\eqref{extended:eqn:vertical:disk}  to obtain
\begin{subequations}\label{FI:vertical:disk}
\begin{align}
(\dot\theta, \dot x, \dot y, \dot \psi) &= \left (\frac{\Pi_\theta}{J}, \frac{\Pi_x\cos\theta}{m}, \frac{\Pi_x\sin\theta}{m}, \frac{\Pi_\psi}{I} \right ),\\
(\dot \Pi_\theta, \dot \Pi_x,\dot \Pi_y,\dot \Pi_\psi)&=( -k_1(\Pi_\theta - \Pi_{\theta 0}), -k_2 (\Pi_x - \Pi_{x0}), -k_3 \Pi_y,-k_4 (\Pi_\psi - \Pi_{\psi 0})). 
\end{align}
\end{subequations}
Then, it is easy to show that for any initial condition the trajectory $\Pi(t)$ will exponentially converge to $(\Pi_{\theta 0}, \Pi_{x0}, 0, \Pi_{\psi 0})$ as $t$ tends to infinity.

\subsection{The Roller Racer}
Consider the roller racer system that is described on pp.42--43 and pp.386--387 in \cite{Bl03}.
The configuration space is $Q = \operatorname{SE}(2) \times \operatorname{S}^1$. Let us use $(\theta, x,y)$ as coordinates for $ \operatorname{SE}(2)$ and $\phi$ for $ \operatorname{S}^1$. Collectively, we use $q = (\theta, x,y,\phi)$ for $Q$ and $p = (p_\theta, p_x,p_y,p_\phi)$ for momentum. The Hamiltonian $H \in C^\infty(T^*Q)$ of the system  is 
\[
H(q,p) = \frac{1}{2m}(p_x^2 + p_y^2) + \frac{(p_\theta - p_\phi)^2}{2I_1} + \frac{p_\phi^2}{2I_2}.
\]
The nonholonomic constraint is
\[
\mathcal C = \left \{\frac{d_1\cos\phi}{I_1} p_\theta  - \frac{ \sin \phi }{m} (p_x\cos\theta + p_y \sin\theta ) + \frac{d_2}{I_2}p_\phi = 0,  -p_x\sin\theta + p_y \cos\theta =0 \right  \}.
\]
Both $H$ and $\mathcal C$ are $\operatorname{SE}(2)$-invariant, so they induce the reduced Hamiltonian $h \in C^\infty( \mathfrak{ se}(2)^* \times T^*\rm{S}^1)$ and the reduced nonholonomic constraint ${\mathcal C}_{\rm red}$ as follows:
\[
h(\Pi,  \phi, p_\phi) = \frac{1}{2m} (\Pi_x^2 + \Pi_y^2) + \frac{(\Pi_\theta - p_\phi)^2}{2I_1} + \frac{p_\phi^2}{2I_2}
\]
and 
\[
{\mathcal C}_{\rm red} = \{ \langle (\Pi,p_\phi), e_1 \rangle = 0,  \langle (\Pi,p_\phi), e_2 \rangle = 0\}
\]
where $\Pi = (\Pi_\theta, \Pi_x, \Pi_y) \in \mathfrak {se}(2)^* \simeq \mathbb R^3$ and the vector fields $e_1$ and $e_2$ on $ \mathbb R^3 \times \operatorname{S^1}$ are given by
\begin{align*}
e_1 &= \frac{1}{\sqrt{\frac{\sin^2\phi}{m}  + \frac{d_1^2\cos^2\phi}{I_1} + \frac{d_2^2}{I_2}}} \left ( \frac{d_1\cos\phi}{I_1}, -\frac{\sin\phi}{m}, 0, -\frac{d_1\cos\phi}{I_1} + \frac{d_2}{I_2}\right ),\\
e_2 &=\frac{1}{\sqrt m}(0,0,1,0).
\end{align*}
Notice that $e_1$ and $e_2$ are orthonormal with respect to the reduced mass tensor of the system given by 
\[
 \begin{bmatrix}
I_1 + I_2 & 0 & 0& I_2 \\
0 & m & 0 & 0 \\
0 & 0 & m & 0 \\
I_2 & 0 & 0 & I_2
\end{bmatrix}.
\]
With respect to this reduced mass tensor, we have
\begin{align*}
e_1^\flat &= \frac{1}{\sqrt{\frac{\sin^2\phi}{m}  + \frac{d_1^2\cos^2\phi}{I_1} + \frac{d_2^2}{I_2}}} ( d_1\cos\phi + d_2, -\sin\phi,0,d_2 ),\\
e_2^\flat &= (0,0,\sqrt m,0).
\end{align*}
Let 
\begin{align*}
J_1(\Pi,\phi, p_\phi) &=  \langle (\Pi,p_\phi),e_1\rangle = \frac{\frac{d_1\cos\phi}{I_1}\Pi_\theta  -\frac{\sin\phi}{m}\Pi_x + (-\frac{d_1\cos\phi}{I_1} + \frac{d_2}{I_2})p_\phi }{\sqrt{\frac{\sin^2\phi}{m}  + \frac{d_1^2\cos^2\phi}{I_1} + \frac{d_2^2}{I_2} }} \\
J_2(\Pi,\phi, p_\phi) &=  \langle (\Pi,p_\phi),e_2\rangle =\frac{\Pi_y}{\sqrt m}.
\end{align*}
Then, the  extended Hamiltonian $\tilde h$ on $\mathfrak {se}(2)^* \times T^*\rm{S}^1$ is computed as
\begin{align*}
\tilde h (\Pi, \phi,p_\phi) &= h (\Pi, \phi,p_\phi) - \frac{1}{2} \sum_{i=1}^2 |J_i(\Pi,p_\phi) |^2\\
&=  \frac{\Pi_x^2}{2m}   + \frac{(\Pi_\theta - p_\phi)^2}{2I_1} + \frac{p_\phi^2}{2I_2}- \frac{\left ( \frac{d_1\cos\phi}{I_1}\Pi_\theta  -\frac{\sin\phi}{m}\Pi_x + (-\frac{d_1\cos\phi}{I_1} + \frac{d_2}{I_2})p_\phi \right )^2}{2\left (\frac{\sin^2\phi}{m}  + \frac{d_1^2\cos^2\phi}{I_1} + \frac{d_2^2}{I_2} \right )}.
\end{align*}
Since $\tilde h$ is independent of $\Pi_y$, the equations of motion of the extended system are written as
\begin{subequations}\label{extended:roller:racer}
\begin{align}
(\dot \theta, \dot x, \dot y) &=  \left (\frac{\partial \tilde h}{\partial \Pi_\theta}, \frac{\partial \tilde h}{\partial \Pi_x} \cos \theta , \frac{\partial \tilde h}{\partial \Pi_x} \sin \theta  \right ), \\
\dot\phi & = \frac{\partial \tilde h}{\partial p_\phi},\\
(\dot \Pi_\theta, \dot \Pi_x, \dot \Pi_y, \dot p_\phi) &= \left ( -\Pi_y \frac{\partial \tilde h}{\partial \Pi_x} , \Pi_y \frac{\partial \tilde h}{\partial \Pi_\theta}, -\Pi_x \frac{\partial \tilde h}{\partial \Pi_\theta}, -\frac{\partial \tilde h}{\partial \phi} \right ) + \sum_{i=1}^2\tilde \lambda_i e_i^\flat,  
\end{align}
\end{subequations}
where
\[
\tilde \lambda_i = \{ \tilde h, J_i\}_{\mathfrak {se}(2)^* \times T^*\rm{S}^1} = \{ \tilde h, J_i\}_{\mathfrak {se}(2)^*} + \{ \tilde h, J_i\}_{T^*\rm{S}^1}
\]
for $i=1,2$. Here $\{\,  , \}_{\mathfrak {se}(2)^*} $ is the minus Lie-Poisson bracket on $\mathfrak{se}(2)^*$ given in \eqref{Poisson:bracket:se2} and $ \{ \, ,\}_{T^*\rm{S}^1}$ is the canonical bracket on $T^*\rm{S}^1$. It is tedious but straightforward to compute the concrete expression of \eqref{extended:roller:racer}, so it is omitted. 

We now construct a feedback integrator for the extended system \eqref{extended:roller:racer}. Choose any number $\tilde h_0$ such that 
\[
\tilde h_0>0
\]
and define a Lyapunov function $V: \mathfrak {se}(2)^* \times T^*\rm{S}^1 \rightarrow \mathbb R$ by
\[
V(\Pi, \phi,p_\phi) = \frac{k_1}{2} |J_1(\Pi,p_\phi)|^2+\frac{k_2}{2} |J_2(\Pi,p_\phi)|^2+ \frac{k_3}{2} |\tilde h(\Pi,\phi,p_\phi) - \tilde h_0|^2
\]
where $k_1, k_2,k_3 >0$. Then, the feedback integrator for \eqref{extended:roller:racer} corresponding to this function $V$ is given as 
\begin{subequations}\label{FI:roller:racer}
\begin{align}
(\dot \theta, \dot x, \dot y) &=  \left (\frac{\partial \tilde h}{\partial \Pi_\theta}, \frac{\partial \tilde h}{\partial \Pi_x} \cos \theta , \frac{\partial \tilde h}{\partial \Pi_x} \sin \theta  \right ),\\
\dot\phi & = \frac{\partial \tilde h}{\partial p_\phi} -\frac{\partial V}{\partial \phi},\label{FI:roller:racer:phi}\\
(\dot \Pi_\theta, \dot \Pi_x, \dot \Pi_y, \dot p_\phi) &= \left ( -\Pi_y \frac{\partial \tilde h}{\partial \Pi_x} , \Pi_y \frac{\partial \tilde h}{\partial \Pi_\theta}, -\Pi_x \frac{\partial \tilde h}{\partial \Pi_\theta}, -\frac{\partial \tilde h}{\partial \phi} \right ) + \sum_{i=1}^2\tilde \lambda_i e_i^\flat - \nabla_{(\Pi,p_\phi)}V, \label{FI:roller:racer:rest}
\end{align}
\end{subequations}
where
 \[
  \nabla_{(\Pi,p_\phi)}V = \left (\frac{\partial V}{\partial \Pi_\theta},\frac{\partial V}{\partial \Pi_x},\frac{\partial V}{\partial \Pi_y},\frac{\partial V}{\partial p_\phi} \right ).
\]
It is  straightforward to compute the partial derivatives of $V$,  so it is omitted.  

\begin{theorem}
There is a positive number $c>0$ such that  every trajectory of \eqref{FI:roller:racer} starting in $\operatorname{SE}(2) \times V^{-1}([0,c])$ remains for all future time in $\operatorname{SE}(2) \times V^{-1}([0,c])$  and asymptotically converges to $\operatorname{SE}(2) \times V^{-1}(0)$ as $t\rightarrow \infty$.
\begin{proof}
It is obvious that $V$ is a first integral of \eqref{FI:roller:racer}. Since $\phi \in \rm{S}^1$ appears in the form of $\cos\phi$ or $\sin\phi$ in the equations of motion, by identifying $\phi$ with $(\cos\phi, \sin\phi)$ it is easy to show that for any $c>0$  the set  $V^{-1}([0,c])$ is a compact subset of $\mathfrak {se}(2)^* \times T^*\rm{S}^1$. 

We now show that the three gradient vectors $\{\nabla J_1, \nabla J_2, \nabla \tilde h\}$ are pointwise linearly independent on $V^{-1}(0)$. Since $\tilde h = h - \frac{1}{2}\sum_{i=1}^2J_i^2$,  the pointwise linear independence of $\{\nabla J_1, \nabla J_2, \nabla  \tilde h\}$ on $V^{-1}(0)$ is equivalent to that of $\{\nabla J_1, \nabla J_2, \nabla  h\}$. It suffices to show the pointwise linear independence of the gradients taken with respect to the momentum variables $(\Pi, p_\phi)$ only, ignoring the configuration variable $\phi$. We have
\[
\begin{bmatrix}
\nabla J_1 \\
\nabla J_2 \\
\nabla h
\end{bmatrix} = \begin{bmatrix}
\frac{d_1\cos\phi}{I_1}& -\frac{\sin\phi}{m} & 0 & -\frac{d_1\cos\phi}{I_1} + \frac{d_2}{I_2} \\
0 & 0 & 1 & 0\\
\frac{\Pi_\theta - p_\phi}{I_1} & \frac{\Pi_x}{m} & \frac{\Pi_y}{m} & \frac{-\Pi_\theta + p_\phi}{I_1} + \frac{p_\phi}{I_2}
\end{bmatrix}
\]
which is transformed through row and column operations to 
\begin{equation}\label{that:matrix}
\begin{bmatrix}
d_1\cos\phi& -\sin\phi& 0 &  d_2 \\
0 & 0 & 1 & 0\\
\Pi_\theta - p_\phi & \Pi_x &0  &  p_\phi 
\end{bmatrix}.
\end{equation}
Suppose that there is a point $(\Pi, \phi, p_\phi) \in V^{-1}(0)$ at which the  matrix in \eqref{that:matrix} has rank less than 3. Then, the determinants of the minor consisting of the last three columns and the minor consisting of the first, third and fourth columns are both zeros, i.e.
\begin{align*}
-p_\phi\sin\phi - d_2 \Pi_x = 0, \quad d_1p_\phi \cos\phi - d_2(\Pi_\theta - p_\phi) =0,
\end{align*}
which are solved for $\Pi_\theta$ and $\Pi_x$ as follows:
\begin{equation}\label{Pithetax:pphi}
\Pi_\theta = \left ( \frac{d_1}{d_2}\cos\phi +1\right )p_\phi, \quad \Pi_x = -\frac{\sin\phi}{d_2}p_\phi,
\end{equation}
Since $(\Pi, \phi, p_\phi) \in V^{-1}(0)$, we have $\Pi_y = J_2(\Pi, \phi, p_\phi) = 0$. Substituting  these in, we get
\[
J_1(\Pi,\phi,p_\phi) = \frac{p_\phi}{d_2}\sqrt{\frac{\sin^2\phi}{m}  + \frac{d_1^2\cos^2\phi}{I_1} + \frac{d_2^2}{I_2} },
\]
which must vanish since $(\Pi, \phi, p_\phi) \in V^{-1}(0)$. Hence,
\[
0 = |J_1(\Pi,\phi,p_\phi) | \geq \frac{1}{\sqrt I_2}|p_\phi|,
\]
which implies $p_\phi = 0$ and thus $\Pi_\theta = \Pi_x = 0$ by  \eqref{Pithetax:pphi}. We now have $(\Pi, \phi, p_\phi) = (0, \phi, 0)$, implying 
\[
V(0, \phi, 0) \geq \frac{k_2}{2}|\tilde h(0,\phi, 0) - \tilde h_0|^2 =   \frac{k_2}{2}|\tilde h_0|^2 >0
\]
which contradicts the assumption that $(\Pi, \phi, p_\phi) \in V^{-1}(0)$. Therefore, the matrix in \eqref{that:matrix} has full rank everywhere on $V^{-1}(0)$, which eventually implies the pointwise linear independence of  $\{ \nabla J_1, \nabla J_2, \nabla \tilde h\}$ on $V^{-1}(0)$.  Then by Theorem 2.5 in \cite{ChJiPe16}, there is a number $c>0$ such that every trajectory of the subsystem \eqref{FI:roller:racer:phi} and \eqref{FI:roller:racer:rest} starting in $V^{-1}([0,c])$ remains for all future time in $V^{-1}([0,c])$  and asymptotically converges to $V^{-1}(0)$ as $t\rightarrow \infty$, from which the theorem follows.
\end{proof}
\end{theorem}

\subsection{The Heisenberg System}
The Hamiltonian of the Heisenberg system is given by
\[
H(q,p)  = \frac{1}{2}|p|^2,
\]
where $q = (x,z,y)\in  \mathbb R^3$ and $p  = (p_x,p_y,p_z)\in \mathbb R^3$. The constraint is given by
\[
0 = p_z - yp_x + xp_y = p\cdot e,
\]
where 
\[
e = (-y,x,1).
\]
We here intentionally did not normalize $e$, but we use the formula \eqref{lambda:general} to compute the multiplier lambda. We compute 
\[
\lambda = 0.
\]
Hence, the equations of motion of the Heisenberg system are given by
\begin{subequations}\label{Heisenberg:original}
\begin{align}
\dot q &= p, \label{Heisenberg:q:original}\\
\dot p &= 0, \label{Heisenberg:p:original}\\
0 &= p\cdot e.
\end{align}
\end{subequations}
Notice that the Hamiltonian $H$ is already a first integral of \eqref{Heisenberg:q:original} and \eqref{Heisenberg:p:original} in the entire phase space due to $\lambda =0$. Hence, there is no need to extend it, but in a sense they are already in an extended form. Let us formally write the extended system without the constraint as follows:
\begin{subequations}\label{Heisenberg:extended}
\begin{align}
\dot q &= p,\\
\dot p &= 0.
\end{align}
\end{subequations}
The extended Hamiltonian is the same as the original Hamiltonian $H$.  Both the extended system and the original system have the following common first integrals: the Hamiltonian $H$, the constraint momentum 
\[
J_1(q,p) = p\cdot e = p_z - yp_x + xp_y,
\]
and the $z$ component of the angular momentum 
\[
J_2(q,p) =   xp_y - yp_x.
\]
Instead of the triple $\{H, J_1,J_2\}$, we could equivalently  use the triple $\{H, p_z, J_2\}$ as a set of first integrals since $p_z = J_1 - J_2$, but we will continue to use the triple $\{H, J_1,J_2\}$. 

Choose any three numbers $H_0$ and $J_{20}$ such that $H_0>0$. Then, define a Lyapunov function $V$ by
\[
V(q,p) = \frac{k_1}{2}|J_1(q,p)|^2 + \frac{k_2}{2}|J_2(q,p) - J_{20}|^2 + \frac{k_3}{2}|H(q,p) - H_0|^2,
\]
where $k_1,k_2,k_3 >0$. It satisfies
\[
V^{-1}(0) = \{ (q,p) \in \mathbb R^3 \times \mathbb R^3 \mid J_1(q,p) = 0, J_2(q,p) = J_{20}, H(q,p) = H_0\},
\]
which is not a compact set. 
The gradient of $V$ is given by
\begin{align*}
\nabla_q V &= k_1 J_1(q,p) (p_y,-p_x,0) + k_2 (J_2(q,p) - J_{20})  (p_y,-p_x,0), \\
\nabla_p V &= k_1 J_1(q,p) (-y,x,1) + k_2 (J_2(q,p) - J_{20}) (-y,x,0) +k_3(H(q,p)-H_0) p.
 \end{align*}
It is easy to see that the function $V$ is a first integral of \eqref{Heisenberg:extended}.
\begin{theorem}
For any $c$ satisfying
\[
0 < c < \min \{ k_3|H_0|^2/2, k_1^2/2k_3\},
\]
the set of critical points of $V$ in $V^{-1}([0,c])$ equals $V^{-1}(0)$.
\begin{proof}
Let $(q,p) \in V^{-1}([0,c])$ be a critical point of $V$. Then it satisfies
\begin{align}
0 &= (k_1 J_1(q,p) + k_2 (J_2(q,p) - J_{20}) ) (p_y,-p_x,0), \label{VqH}\\
0&= k_1 J_1(q,p) (-y,x,1) + k_2 (J_2(q,p) - J_{20}) (-y,x,0) +k_3(H(q,p)-H_0) p.\label{VpH}
 \end{align}
 First, consider the case where $p_x \neq0$. Then, \eqref{VqH} implies $k_1 J_1(q,p) + k_2 (J_2(q,p) - J_{20}) =0$. Substitute this in \eqref{VpH}, and then we obtain $ (H(q,p) - H_0) p_x = 0$ and $k_1J_1(q,p) + k_3 (H(q,p) - H_0)p_z =0$. Since $p_x \neq0$, it follows that $H(q,p) - H_0=0$, and thus $J_1(q,p) =0$. Hence, $J_2(q,p)  - J_{20} = 0$. It implies that $(q,p) \in V^{-1}(0)$. The case of $p_y \neq0$ similarly leads to $(q,p) \in V^{-1}(0)$.  Now consider that case where $p_x = p_y = 0$. If $p_z \neq0$, then, the third component of the vector equation \eqref{VpH} implies $k_1 + k_3 (H(q,p) - H_0) = 0$, so $H(q,p) - H_0 = -k_1/k_3$. Then, $V(q,p) \geq k_3|H(q,0) - H_0|^2 /2\geq k_1^2/2k_3 >c$, contradicting $(q,p) \in V^{-1}([0,c])$.   If $p_z = 0$, then $V(q,p) \geq k_3|H(q,0) - H_0|^2/2 \geq k_3|H_0|^2/2 >c$, contradicting $(q,p) \in V^{-1}(0)$. Thus, the case $p_x =p_y = 0$ is not possible. Considering the above arguments, we come to the conclusion that every critical point of $V$ in $V^{-1}([0,c])$ is contained in $V^{-1}(0)$. Since $0$ is the minimum value of $V$, every point of $V^{-1}(0)$ is a critical point of $V$. Therefore, the set of critical points of $V$ in $V^{-1}([0,c])$ equals $V^{-1}(0)$.
  \end{proof}
\end{theorem}

The feedback integrator corresponding to $V$ is given by
\begin{subequations}
\begin{align*}
\dot q &= p - \nabla_qV,\\
\dot p &=  - \nabla_p V.
\end{align*}
\end{subequations}
Although we do not have a convergence proof  due to non-compactness of $V^{-1}(0)$, we expect that it will perform well practically as was the case for the knife edge system. Numerical simulation is left to the reader.

\subsection{The Nonholonomic Oscillator}
We design a feedback integrator and a Lagrange-d'Alembert integrator for the nonholonomic oscillator and compare the performances of the two integrators.  We refer the reader to \cite{MoVeSwede} on the  properties, including integrability, of the nonholonomic oscillator.

\paragraph{Feedback Integrator.}
The Hamiltonian of the nonholonomic oscillator is given by
\[
H = \frac{1}{2}|p|^2 + \frac{1}{2}|q|^2
\]
where $q = (x,y,z) \in \mathbb R^3$ and $p = (p_x, p_y, p_z) \in \mathbb R^3$, and the nonholonomic constraint set is given by
\[
{\mathcal C} = \{ p_x + y p_z = 0\} = \{ p\cdot e = 0\}
\]
where 
\[
e = \frac{1}{\sqrt{1+y^2}}(1,0,y). 
\]
The equations of motion of the system are
\begin{subequations}
\begin{align*}
\dot q &= p ,\\
\dot p &=  - q + \lambda e,\\
0&=p\cdot e.
\end{align*}
\end{subequations}
Let
\[
J(q,p) = p \cdot e = \frac{p_x + y p_z}{\sqrt{1+y^2}}.
\]
The extended Hamiltonian $\tilde H$ is computed as
\[
\tilde H (q,p) = \frac{(yp_x - p_z)^2}{2(1+y^2)} + \frac{1}{2}p_y^2 + \frac{1}{2}|q|^2,
\]
and the equations of motion of the extended system are given by
\begin{subequations}
\begin{align*}
\dot q &= \nabla_p \tilde H,\\
\dot p &= -\nabla_q \tilde H + \tilde \lambda e,
\end{align*}
\end{subequations}
where
\begin{align*}
\nabla_q \tilde H &= \left  (x, y + \frac{(yp_x - p_z)(p_x + yp_z)}{(1+y^2)^2}, z \right ) \\
\nabla_p \tilde H&= \left (\frac{(yp_x - p_z)y}{1+y^2}, p_y, \frac{p_z - yp_x}{1+y^2} \right ),
\end{align*}
and
\[
\tilde \lambda = \frac{(x + yz)(1+y^2) + (p_x y - p_z)p_y}{(1+y^2)^{3/2}}.
\]
This system has the following three constants of motion: the Hamiltonian $\tilde H$, and  the constraint momentum $J(q,p)$, and the energy of the $y$ dynamics defined by
\[
 H_y(q,p) = \frac{1}{2}p_y^2 + \frac{1}{2}y^2.
\]
Take any two numbers $\tilde H_0$ and $H_{y0}$  so that $0< H_{y0} \leq \tilde H_0$, and let
\[
V(q,p) = \frac{k_1}{2} |\tilde H(q,p) - \tilde H_0|^2 + \frac{k_2}{2}|H_{y0}(q,p) - H_{y0}|^2 + \frac{k_3}{2}|J(q,p)|^2,
\]
where $k_1$, $k_2$ and $k_3$ are positive constants. Then,  the feedback integrator corresponding to this function $V$ is given by
\begin{subequations}
\begin{align*}
\dot q &= \nabla_p \tilde H - \nabla_qV,\\
\dot p &= -\nabla_q \tilde H + \tilde \lambda e - \nabla_pV,
\end{align*}
\end{subequations}
where
\begin{align*}
\nabla_qV &= k_1 \Delta \tilde H \nabla_q \tilde H + k_2 \Delta H_y \nabla_qH_y + k_3 J \nabla_q J\\
\nabla_p V &= k_1 \Delta \tilde H \nabla_p \tilde H + k_2 \Delta H_y \nabla_pH_y + k_3  J \nabla_p J
\end{align*}
and
\begin{align*}
\nabla_q H_y &= (0,y,0), \quad \nabla_p H_y = (0,p_y,0),\\
\nabla_q J &= \left (0, \frac{p_z-yp_x}{(1+y^2)^{3/2}},0 \right ), \quad \nabla_p J =  \left ( \frac{1}{\sqrt{1+y^2}}, 0, \frac{y}{\sqrt{1+y^2}} \right )
\end{align*}
with $\Delta \tilde H = \tilde H(q,p) - \tilde H_0$ and $\Delta H_y = H_y(q,p) - H_{y0}$.
 
\paragraph{Lagrange-d'Alembert Integrator.} Consider a  discrete Lagrange-d'Alembert (DLA) integrator for the nonholonomic oscillator (\cite{McPe}). Here we just take the DLA in the time step $h$. With the map $\phi: \mathbb{R}^{3}\times \mathbb{R}^{3}\rightarrow T\mathbb{R}^{3}$ given by
\begin{equation}
\label{EZPHI}
\phi(q_{0},q_{1})=\left (q_{0},\frac{q_{1}-q_{0}}{h} \right )
\end{equation}
and letting $L_{d}=L\circ \phi$, we have
\begin{equation}
\label{DiscLag}
L_{d}(q_{0},q_{1})=\frac{1}{2}\left |\frac{q_{1}-q_{0}}{h}\right |^{2}-\frac{1}{2}|q_{0}|^{2}
\end{equation}
so that
\begin{align}
D_{1}L_{d}(q_{0},q_{1})&=\frac{q_{0}-q_{1}}{h^{2}}-q_{0} \\
D_{2}L_{d}(q_{0},q_{1})&=\frac{q_{1}-q_{0}}{h^{2}}
\end{align}
The DLA equations, recall, are, for $i\in\{0,\dots N-2\}$,
\begin{align}
D_{2}L_{d}(q_{i},q_{i+1})+D_{1}L_{d}(q_{i+1},q_{i+2})&=\lambda_{i}A(q_{i+1})\\
(q_{i+1},q_{i+2})&\in C_{d}
\end{align}
where $C_{d} \subset Q\times Q$ is the discrete constraint distribution. For our choice of $\phi$ \eqref{EZPHI} this gives
\begin{equation}
C_{d}=\{(q_{0},q_{1}) \, | \, A(q_{0})\cdot (q_{1}-q_{0})=0\}
\end{equation}
which for $A=dx+ydz$ or $A = (1,0,y)$, gives
\begin{equation}
C_{d}=\{(q_{0},q_{1}) \, | \, x_{1}-x_{0}+y_{0}(z_{1}-z_{0})=0\}.
\end{equation}
This all translates to the following set of equations:
\begin{align}
\label{ezDLA}
\frac{q_{i+1}-q_{i}}{h^{2}}+\frac{q_{i+1}-q_{i+2}}{h^{2}}-q_{i+1}&=\lambda \begin{bmatrix} 1\\ 0\\ y_{i+1}\end{bmatrix},\\
x_{i+2}-x_{i+1}+y_{i+1}(z_{i+2}-z_{i+1})&=0.
\end{align}
We can solve these equations explicitly by first isolating $q_{i+2}$, and then plugging this expression into the second equation above to determine $\lambda$.
This procedure yields,
\begin{equation}
\label{dynamics}
q_{i+2}=2q_{i+1}-q_{i}-h^{2}q_{i+1}-h^{2}\lambda \begin{bmatrix} 1\\ 0\\ y_{i+1}\end{bmatrix}
\end{equation}
which plugged into the constraint condition and isolating $\lambda$ gives
\begin{equation}
\label{EzLambda}
\lambda=\frac{(1-h^{2})x_{i+1}-x_{i}+y_{i+1}((1-h^{2})z_{i+1}-z_{i})}{h^{2}(1+y_{i+1}^{2})}.
\end{equation}
These two equations \eqref{dynamics} and \eqref{EzLambda} then explicitly give the discrete dynamics.

\begin{figure}[t]
\hspace{-15mm}
\includegraphics[scale = 0.35]{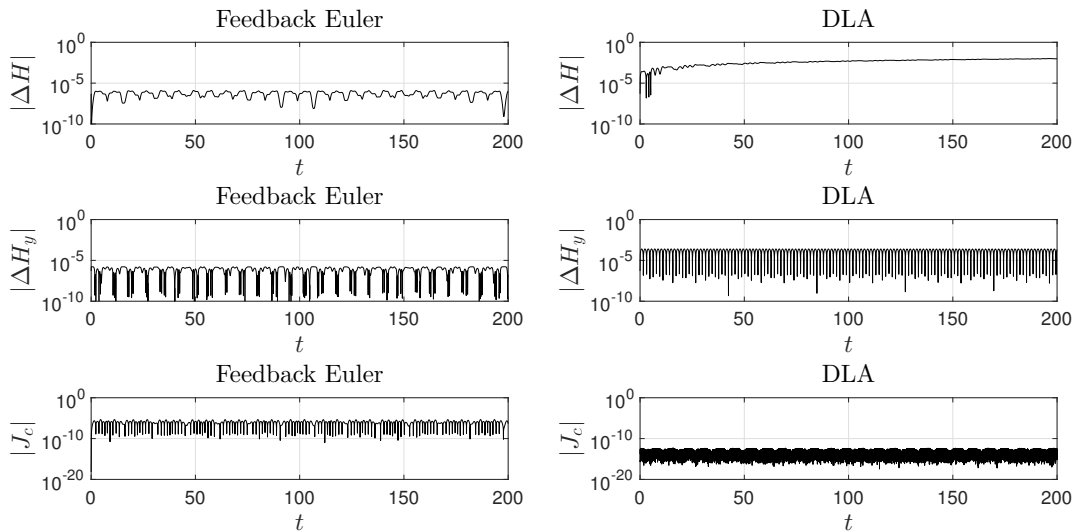}
\caption{The left column shows the total energy error $\Delta H$, the $y$-energy error $\Delta H_y$ and the nonholonomic constraint $J_c(q,p) =  p_x + yp_z$ of the numerical solution  of the nonholonomic oscillator system generated by  a feedback integrator with the Euler scheme with step size $\Delta t = 10^{-3}$ over the time interval $[0,200]$. The right column shows those generated by a first-order DLA integrator with step size $h = 10^{-3}$ over the time interval $[0,200]$.}
\label{figure:oscillator}
\end{figure}

\paragraph{Comparisons.} Take time step size $\Delta t =h =  10^{-3}$ for both integrators. Take the  initial state
\begin{equation}\label{ns:ic}
q(0) = (1,1,1), \quad p(0,0,0) = (0,0,0)
\end{equation}
for the feedback integrator, which  translates to 
\begin{equation}\label{ns:dla:ic}
q_0 = (1,1,1), \quad q_1 = (1,1,1)
\end{equation}
for the DLA integrator. Choose the following values for $k_1$, $k_2$ and $k_3$ of the feedback integrator:
\begin{equation}\label{no:gains}
k_1 = 500, \quad k_2 = 300, \quad k_3 = 300.
\end{equation}
The usual Euler scheme is used for the feedback integrator. The simulations are run over the time interval $[0, 200]$.  The discrete momentum is computed as
\[
p_{i+1} = (q_{i+1} - q_i)/h.
\]
The total energy error $\Delta H = H(q,p) - H(q(0),p(0))$, the $y$-energy error, $\Delta H_y = H_y(q,p) - H_y(q(0),p(0))$ and the nonholonomic constraint $J_c (q,p) = p_x + yp_z$ are plotted in Figure \ref{figure:oscillator}. In the figure,  it is observed that the DLA integrator is poor at preserving the energy while it preserves the other two conserved quantities well. In contrast, the feedback  integrator preserves all the three conserved quantities well. Notice that both of the schemes under comparison are  first-order integration schemes.

\begin{figure}[t]
\hspace{-15mm}
\includegraphics[scale = 0.29]{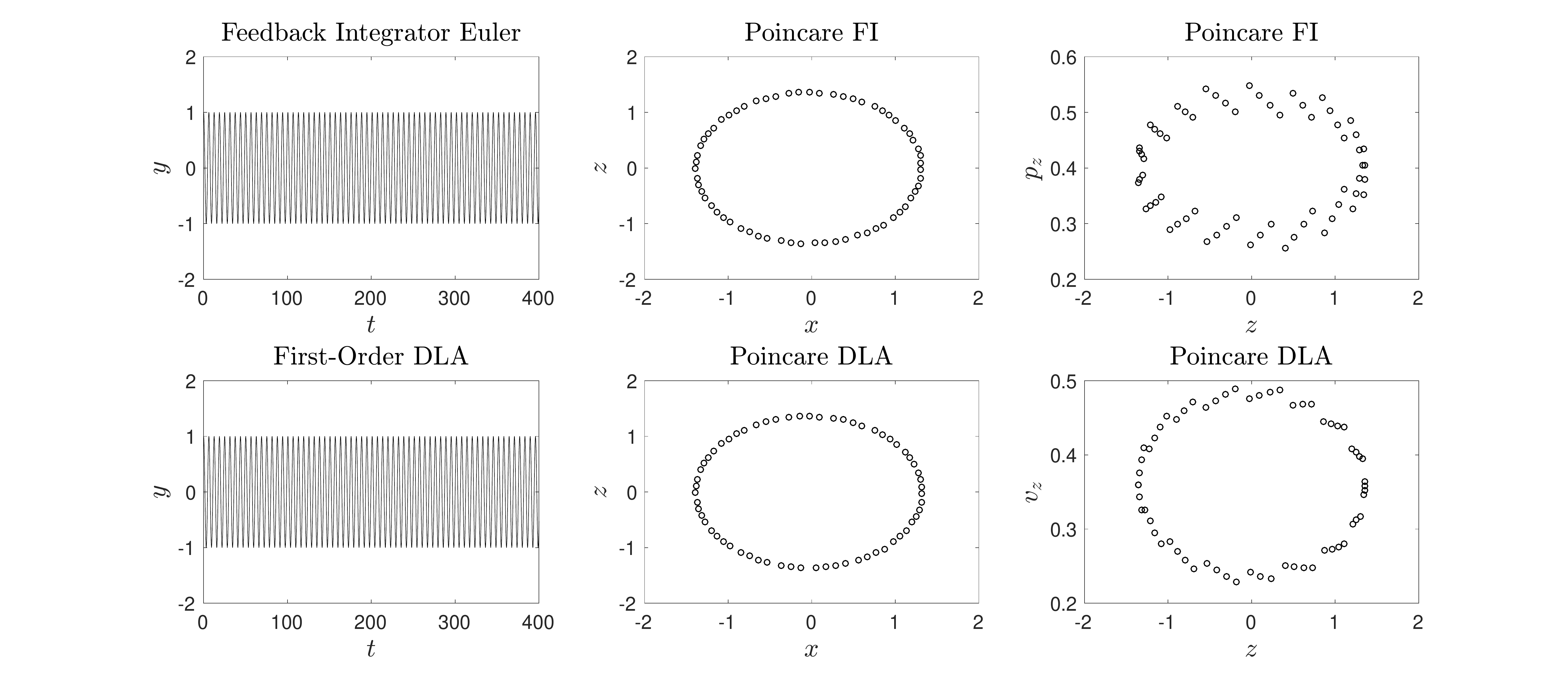}
\caption{The first row plots  $y(t)$ of the nonholonomic oscillator trajectory, the snapshots of the trajectory in the $xz$ plane and the snapshots  of the trajectory in the $zp_z$ plane when $y(t)$ crosses zero with $p_y(t)>0$ for the feedback integrator with the Euler scheme.  The second row plots  $y(t)$ of the nonholonomic oscillator  trajectory, the snapshots of the trajectory in the $xz$ plane and the snapshots  of the trajectory in the $zv_z$ plane when $y(t)$ crosses zero with $v_y(t)>0$ for  the first-order DLA method. }
\label{figure:Poincare:FI}
\end{figure}

We now compare the two methods by plotting Poincar\'e maps; refer to \cite{MoVeSwede}  on the integrability of the nonholonomic oscillator.  This time we run the simulations over the time interval $[0,400]$ with the other conditions being the same as in \eqref{ns:ic} -- \eqref{no:gains}.   The $y$-dynamics of the nonholonomic oscillator are 
\[
\dot y = p_y, \quad \dot p_y = -y.
\]
With the initial condition \eqref{ns:ic}, the $y$ trajectory becomes $y(t) = \cos t$ which is periodic with period $2\pi$. We take two snapshots of the trajectory each time $y(t)$ crosses zero with $\dot y(t)>0$: one snapshot in the $xz$-plane and the other in the $zp_z$-plane (or, $zv_z$-plane). The results are plotted in \ref{figure:Poincare:FI}. In the first row in the figure, plotted are the $y(t)$ trajectory, the snapshots  in the $xz$-plane and the snapshots  in the $zp_z$-plane for the feedback integrator method.  In the second row, plotted are 
the $y(t)$ trajectory, the snapshots  in the $xz$-plane and the snapshots  in the $zv_z$-plane for the first-order DLA method. It is observed that both methods produce similar results.

\section{Conclusion}
We have successfully developed a theory of feedback integrators for nonholonomic mechanical systems with or without symmetry, where the case with symmetry was studied in the cases that configuration space is a symmetry group or a trivial principal bundle. We have successfully applied the nonholonomic integrators to the following systems: the Suslov problem on $\operatorname{SO}(3)$, the knife edge, the Chaplygin sleigh, the vertical rolling disk, the roller racer,  the Heisenberg system, and the nonholonomic oscillator.  We plan to develop feedback integrators  for nonholonomic mechanical systems with symmetry on non-trivial principal bundles and mechanical systems with holonomic constraints in a future work.


\end{document}